\documentclass[10pt]{article}

\usepackage{amsfonts} 
\usepackage{amsmath,hhline,latexsym,mathrsfs}
\usepackage{amssymb}
\usepackage{relsize}
\DeclareMathAlphabet\mathbfcal{OMS}{cmsy}{b}{n}

\usepackage[latin1]{inputenc}

\usepackage{hyperref}

\usepackage{color}
\usepackage{graphicx}
\usepackage{comment}

\newcommand{\bfin}{\hspace{12cm}\ \ \ \ \rule[-10pt]{6pt}{6pt}}

\newtheorem{theorem}{\sc Theorem}[section]
\newtheorem{lemma}{\sc Lemma}[section]
\newtheorem{definition}{\sc Definition}[section]
\newtheorem{proposition}{\sc Proposition}[section]

\newtheorem{remark}{Remark}

\newcommand{\eps}{\varepsilon}

\newcommand{\ph}{\varphi}

\newcommand{\N}{\mbox{$I \kern -4pt N$}}
\newcommand{\Q}{\mbox{$Q \kern -8pt I$}}
\newcommand{\R}{\mbox{$I \kern -4pt R$}}
\newcommand{\C}{\mbox{$C \kern -8pt I$}}

\def \beq {\begin{equation}}
\def \eeq {\end{equation}}
\def \ba {\begin{array}}
\def \ea {\end{array}}

\providecommand{\tabularnewline}{\\}

\def\R{{\bf R}}
\usepackage[english,ruled,vlined]{algorithm2e}

\title{Singular asymptotic expansion of the exact control for a linear model of the Rayleigh beam}
\author{Carlos \textsc{Castro} \thanks{Departamento de Matem\'atica e Inform\'atica Aplicadas a la Ingeniería Civil y Naval, Universidad Polit\'ecnica de Madrid, Madrid, 28040, Spain. E-mail: {\tt carlos.castro@upm.es}}
\and 
Arnaud \textsc{M\"unch}\thanks{Laboratoire de Math\'ematiques Blaise Pascal, Universit\'e Clermont Auvergne, UMR CNRS 6620, Campus universitaire des C\'ezeaux, 3, place Vasarely, 63178, Aubi\`ere, France. E-mail: {\tt arnaud.munch@uca.fr.}} 
}

\begin{document}

\maketitle

\begin{abstract}
The Petrowsky type equation $y_{tt}^\eps+\eps y_{xxxx}^\eps - y_{xx}^\eps=0$, $\eps>0$ encountered in linear beams theory is null controllable through Neumann boundary controls. Due to the boundary layer of size of order $\sqrt{\eps}$ occurring at the extremities, these boundary controls get singular as $\eps$ goes to $0$. Using the matched asymptotic method, we describe the boundary layer of the solution $y^\eps$ then derive a rigorous second order asymptotic expansion of the control of minimal $L^2-$norm, with respect to the parameter $\eps$. In particular, we recover that the leading term of the expansion is a null Dirichlet control for the limit hyperbolic wave equation, in agreement with earlier results due to J-.L. Lions in the eighties. Numerical experiments support the analysis.     
\end{abstract}
\vskip .5cm\noindent
{\bf Keywords~:} Singular controllability, Boundary layers, Asymptotic analysis, Numerical experiments.

\par\noindent
{\bf Mathematics Subject Classification~:} 93B05, 58K55.


\section{Introduction}

Let $\Omega=(0,1)$, $\Gamma=\{1\}$. For any $T>0$, we note $Q_T:=\Omega\times (0,T)$, 
$\Gamma_T:=\Gamma\times (0,T)$ and $\Sigma_T:=\partial\Omega\times (0,T)$. This work is concerned with the null controllability property with respect to the parameter $\eps>0$ of the following linear equation of Petrowsky type
\begin{equation}\label{eq:sys_yeps}
\left\{
\begin{aligned}
& y^{\eps}_{tt}+ \eps y_{xxxx}^{\eps} - y_{xx}^{\eps}=0, &\text{ in }  \, Q_T, \\
&  y^{\eps}(0,\cdot)=y^{\eps}(1,\cdot)=y_x^{\eps}(0,\cdot)=0,  \quad y_x^{\eps}(1,\cdot)=v^{\eps} &\text{ in }  \, (0,T),   \\ 
&  (y^{\eps}(\cdot,0),y_t^{\eps}(\cdot,0))=(y_0,y_1), &\text{ in }  \, \Omega.
\end{aligned}
\right.
\end{equation}
Here, $v^{\eps}$ is a control function in $L^2(0,T)$. This system is a simplified version of the modelization of the dynamic of a linear isotropic beam of length and mass equal to one. $\eps$ stands for the product of the young modulus with the area moment inertial of the cross section (we refer to \cite{Miklowitz}). $y^{\eps}=y^{\eps}(x,t)$ is the deflection of the beam at point $x\in \Omega$ and time $t\in (0,T)$. $y_0$ denotes the initial position and $y_1$ the initial velocity assumed in $L^2(\Omega)$ and $H^{-2}(\Omega)$ respectively.

For any $\eps>0$, $v^{\eps}\in L^2(\Omega)$ and $(y_0,y_1)\in L^2(\Omega)\times H^{-2}(\Omega)$, there exists exactly one solution $y^{\eps}$ to (\ref{eq:sys_yeps}), with the regularity $y^{\eps}\in C^0([0,T],L^2(\Omega))\cap C^1([0,T], H^{-2}(\Omega))$ (see \cite{lions2}, chapter 4) and the following estimate: 
\begin{equation} 
\nonumber
\Vert y^{\eps}\Vert_{L^{\infty}(0,T;L^2(\Omega))}+\Vert y_t^{\eps}\Vert_{L^{\infty}(0,T;H^{-2}(\Omega))}\leq c_{\eps}\biggl(\Vert y_0\Vert_{L^2(\Omega)}+\Vert y_1\Vert_{H^{-2}(\Omega)}+\Vert v^{\eps}\Vert_{L^2(0,T)} \biggr)
\end{equation}
for some constant $c_{\eps}>0$. 

Accordingly, for any final time $T>0$, the associated null controllability problem at time $T$ is the following: for each $(y_0,y_1)\in L^2(\Omega)\times H^{-2}(\Omega)$, find a control function $v^{\eps}\in L^2(0,T)$ such that the corresponding solution to (\ref{eq:sys_yeps}) satisfies 
\begin{equation}
(y^{\eps}(\cdot,T),y_t^{\eps}(\cdot,T))=(0,0) \;\; \mbox{ in } L^2(\Omega)\times H^{-2}(\Omega).  \label{eq:null}
\end{equation}
For any $\eps>0$, existence of null controls is proved in \cite{lions2} assuming that $T>2$. Precisely, the following observability inequality is obtained for some constant $C>0$ independent of $\eps$
\begin{equation} \label{eq:obsineqeps}
\Vert \ph^\eps_0\Vert^2_{H^1_0(\Omega)}+\Vert \ph^\eps_1\Vert^2_{L^2(\Omega)}+ \eps \Vert \ph_{0,xx}^\eps\Vert^2_{L^2(\Omega)}\leq C\int_0^T \eps \vert \ph_{xx}^\eps(1,\cdot)\vert^2,\, \forall (\ph^\eps_0,\ph^\eps_1)\in H^2_0(\Omega)\times L^2(\Omega)
\end{equation}
where $\ph^{\eps}$ solves the corresponding homogeneous adjoint associated to the initial condition ($\ph^\eps_0,\ph^\eps_1$), see \eqref{eq:sys_pheps1d}. Let us introduce the non-empty set of null controls
\begin{equation}\nonumber
\mathcal{C}(y_0,y_1,T,\eps):=\{(y,v): v\in L^2(0,T); y\,\, \textrm{solves}\,\, (\ref{eq:sys_yeps})\,\, \textrm{and satisfies}\,\, (\ref{eq:null}) \}.
\end{equation}

Since the physical parameter $\eps$ is small with respect to one, the issue of the asymptotic behavior of elements of $\mathcal{C}$ as $\eps$ is smaller and smaller arises naturally. It turns out that the system (\ref{eq:sys_yeps}) is not uniformly controllable with respect to $\eps$ in the sense that the control of minimal $L^2(0,T)$-norm is not uniformly bounded. The following result is proved in \cite[chapter 6]{lions2}, assuming additional regularity on the initial velocity. We also refer to \cite{LionsPerturbation87}.
\begin{theorem}[Lions \cite{lions2}]
\label{theoremLions}
Assume that the initial condition $(y_0,y_1)$ belongs to $L^2(\Omega)\times H^{-1}(\Omega)$ and $T\geq 2$. For any $\eps>0$, let $v^{\eps}$ be the control of minimal $L^2(0,T)$-norm for $y^{\eps}$ solution of (\ref{eq:sys_yeps}). Then, one has 
\begin{equation}
\label{con_vep}
\begin{aligned}
& -\sqrt{\eps}v^\eps \to v \quad \textrm{in}\quad  L^2(0,T)-\textrm{weak}, \quad \textrm{as}\quad \eps\to 0, \\
& y^{\eps} \to y \quad \textrm{in}\quad L^{\infty}(0,T; L^2(\Omega))-\textrm{weak-star} , \quad \textrm{as}\quad \eps\to 0
\end{aligned}
\end{equation}
where $v$ is the control of minimal $L^2(0,T)$-norm for $y$, solution in  $C^0([0,T]; L^2(\Omega))\times  C^1([0,T];H^{-1}(\Omega))$ of the following system :
\begin{equation}\label{eq:sys_y}
\left\{
\begin{aligned}
& y_{tt}- y_{xx}=0, &\text{ in }  \, Q_T, \\
& y(0,\cdot)=0, \, y(1,\cdot)=v,&\text{ in }  \, (0,T) ,   \\
&  (y(\cdot,0),y_t(\cdot,0))=(y_0,y_1), &\text{ in }  \, \Omega.
 \end{aligned}
\right.
\end{equation}
\end{theorem}

Theorem \ref{theoremLions} is obtained by using \textit{a priori} estimates on the solutions and by taking the weak limit in the optimality system characterizing the control of minimal $L^2$-norm. {\color{black} This argument holds in higher dimensions too but here we focus on the one-dimensional case for which we can give a detailed asymptotic expansion of the controls.}

As mentioned in \cite{lions2}, this controllability problem is singular in the following two meanings:  firstly, the convergence result (\ref{con_vep}) holds for $\sqrt{\eps}v^\eps$ and not for $v^\eps$; secondly,  the Neumann null type control $v^{\eps}$ for (\ref{eq:sys_yeps}) degenerates as $\eps\to 0$ into a Dirichlet type control $v$ for the limit problem (\ref{eq:sys_y}). 
The control of minimal norm being of the form $v^{\eps}=\ph_{xx}^{\eps}(1,\cdot)$ on $(0,T)$, this singularity is related to the boundary layer of length $\sqrt{\eps}$ which occurs on the adjoint solution $\ph^{\eps}$ in the neighborhood of the point $x=1$ as $\eps$ tends to $0$.

On the other hand, it is interesting to notice that the spectrum of the underlying operator satisfies a uniform gap property with respect to the parameter $\eps$ \cite{Rellich}. {\color{black} Therefore, inequality (\ref{eq:obsineqeps}) can be also obtained using the classical approach based on reducing the observability inequality to a moment problem for the associated family of exponentials.  The $\varepsilon$ term in the right hand side of (\ref{eq:obsineqeps}) comes from the boundary layer that affects to the boundary observability for the eigenfunctions.}  

Remark also that, without compatibility conditions at $\partial\Omega\times \{0\}$ between the initial and boundary conditions, the solution of (\ref{eq:sys_yeps}) develops as $\eps\to 0$ additional internal boundary layers along the characteristics of (\ref{eq:sys_y}), namely $C_l=\{(x,t)\in Q_T, x-t=0\}$ and $C_d=\{(x,t)\in Q_T, x+t-1=0\}$. {\color{black}  This case would require an specific asymptotic analysis which capture this behavior and we do not consider it here.  Thus we will add compatibility conditions to the initial data when required to avoid this situation}

The main reason of this work is to perform an asymptotic analysis of the controls with respect to the parameter $\eps$ in order notably to precise the rate of convergence of the pair $(\sqrt{\eps}v^{\eps},y^{\eps})$. Such analysis is actually mentioned in the open problems section in \cite[chapter 7]{lions2}. Precisely, we want to characterize the terms $v^k$, $k\geq 0$, in the expansion 
\begin{equation}
\eps^{1/2}v^{\eps}= v^0+\eps^{1/2}v^1+ \eps v^2+\cdots. \label{expansion_control}
\end{equation}

For any fixed function $v^{\eps}$, the asymptotic analysis of solutions of singular systems like (\ref{eq:sys_yeps}) can be performed using the matched asymptotic method \cite{ColeBook96, eckhausbook79}. Precisely, under additional regularity and compatibility assumptions on the data $(y_0,y_1)$ and on the function $v^{\eps}$, one may construct explicitly strong convergent approximations of $y^{\eps}$.
We refer for instance to \cite{AmiratMunch} for an advection-diffusion equation. When the exact controllability issue comes into play, such analysis requires more care, since the regularity/compatibility conditions mentioned above become additional constraints on the control set $\mathcal{C}$. Typically, $L^2$ regularity for the control is in general not sufficient to get strong convergent results for the corresponding controlled solution, at any order. It is then necessary to enrich the set $\mathcal{C}$ and to modify the optimality system accordingly. This is probably the reason for which there exists in the literature only very few asymptotic analysis for controllability problem, \textit{a fortiori} for singular partial differential equations. We mention \cite{Marbach}, \cite{AmiratMunch} following \cite{CoronGuerrero}.
We also mention the book \cite{Lions1973} and the review \cite{kurina} in the close context of optimal control problems. 

Hopefully, as noticed in \cite{Lebeau, ErvedozaZuazua}, a small modification of the set $\mathcal{C}$, still in an $L^2$- framework, allows to recover \textit{a posteriori} smooth controls for smooth data. This property, which do not hold in the context of the advection-diffusion equation discussed in  \cite{AmiratMunch, CoronGuerrero}, is a crucial point for the analysis below. 
This modification also allows to impose simply appropriate compatibility conditions on $\partial\Omega\times (0,T)$ and therefore avoid the internal layers mentioned above along the characteristics $C_l$ and $C_d$. 

Our analysis relies on the matched asymptotic expansion which allows to construct explicitly an expansion of $y^{\eps}$ and $v^{\eps}$. 
The control of minimal $L^2$-norm $v^{\eps}$ for (\ref{eq:sys_yeps}) (precisely, in view of Theorem \ref{theoremLions}, the control which minimizes $v\to \Vert \sqrt{\eps} v^{\eps} \Vert_{L^2(0,T)}$ over $\mathcal{C}(y_0;y_1,\eps,T)$) is given by $v^{\eps}:=\ph^{\eps}_{xx}$ where $\ph^{\eps}$ solves the homogeneous problem
\begin{equation}\label{eq:sys_pheps1d}
\left\{
\begin{aligned}
&\ph^{\eps}_{tt}+ \eps \, \ph^{\eps}_{xxxx} -  \ph^{\eps}_{xx}=0, &\text{ in }  \, Q_T, \\
&\ph^{\eps}=\ph_x^{\eps}=0,&\text{ on }  \,\Sigma_T,   \\ 
&(\ph^{\eps}(\cdot,0),\ph_t^{\eps}(\cdot,0))=(\ph^{\eps}_0, \ph^{\eps}_1),&\text{ in }  \, \Omega,
\end{aligned}
\right.
\end{equation}
associated to the initial condition $(\ph^{\eps}_0, \ph^{\eps}_1)$ minimizer of the so-called conjugate functional $J_{\eps}^{\star}:H^2_0(\Omega)\times L^2(\Omega)\to \mathbb{R}$ defined by 
\begin{equation}\nonumber
J_{\eps}^{\star}(\ph^{\eps}_0,\ph^{\eps}_1)=\frac{\eps}{2}\Vert \ph_{xx}^{\eps}(1,\cdot)\Vert^2_{L^2(0,T)}- (y_0,\ph_1^{\eps})_{L^2(\Omega),L^2(\Omega)}  + (y_1,\ph_0^{\eps})_{H^{-2}(\Omega),H^2(\Omega)}. 
\end{equation}
We clearly see here that the singular character of the control $v^{\eps}$ is due to the boundary layer occurring at $x=1$ on the quantity $\ph_{xx}^\eps(1,\cdot)$. In particular, such singular character does not occur for a distributed control in $\Omega$.  
%

We also emphasize that the asymptotic analysis with respect to $\eps$ is also relevant from an approximation viewpoint. Precisely, an approximation of $v^{\eps}$ is usually obtained through the minimization of the conjugate functional (see for instance \cite{munch_arch2010} in a similar context). However as $\eps$ goes to zero, the norm $\Vert \sqrt{\eps}\ph_{xx}^{\eps}\Vert_{L^2(0,T)}$ is not anymore equivalent to $\Vert \ph_0^\eps,\ph_1^\eps\Vert_{H_0^2(\Omega)\times L^2(\Omega)}$ so that the minimization of $J^\star_{\eps}$ becomes ill-conditionned. On the contrary, the characterization of each terms in the expansion (\ref{expansion_control}) leads to well-conditionned extremal problems (at the successive orders of $\eps$). This guarantees an accurate approximation of the control of minimal $L^2$-norm $v^{\eps}$ with respect to the approximation parameter, for arbitrarily small value of $\eps$. {\color{black} This is illustrated below with some numerical experiments.}

This document is organized as follows. In Section \ref{weightedcontrol}, we set the problem in a suitable functional framework and check that smooth controls can be achieved for smooth data, both for (\ref{eq:sys_yeps}) and the limit wave equation. In Section \ref{asymptotic_expansion}, we use the matched asymptotic method and determine expansion of solutions of (\ref{eq:sys_yeps})  and the associated adjoint system (\ref{eq:sys_pheps1d}). In Section \ref{optimality_system}, we substitute these expansion in the optimality system associated to the control of minimal $L^2$-norm and fully characterize each term in the expansion (\ref{expansion_control}). Section \ref{convergence_control} justifies these expansion with convergence results with respect to $\eps$. Eventually, we discuss some experiments in Section \ref{experiments} and provides some perspectives in the concluding section.

\section{Minimal $L^2$-weighted controls}\label{weightedcontrol}

In this section, {\color{black} following \cite{ErvedozaZuazua},} we introduce a class of controls for both system (\ref{eq:sys_yeps}) and the limit one (\ref{eq:sys_y}) which are smooth for smooth data. Note that, in general, this is not the case for controls with minimal $L^2$-norm and the asymptotic analysis below cannot be justified for such controls. 

We introduce a weight function with suitable properties that will be used below to define a class of smooth controls. For $T>2$, the minimal time to have controllability of system (\ref{eq:sys_yeps}), we consider a positive smooth weight function $\eta \geq 0$ with compact support in $(0,T)$, i.e. $\eta \in C_0^\infty (0,T)$, and such that $\eta(t)>\eta_0>0$ in a subinterval $[\delta,T-\delta]\subset (0,T)$ with $\delta $ such that $T-2\delta >2$.

We divide the rest of this section in two subsections where we consider separately the results for systems (\ref{eq:sys_yeps}) and (\ref{eq:sys_y}) respectively.

\subsection{Beam system}

Let $X=L^2(0,1)\times H^{-2}(0,1)$ and $X^*=H^2_0(0,1)\times L^2(0,1)$ its dual, with duality product given by 
\begin{equation} \label{eq_dua_pr}
<(y_0,y_1),(\varphi_0,\varphi_1)>_{X,X^*} = \int_\Omega y_0\varphi_1 \; dx - (y_1,\varphi_0)_{H^{-2},H^{2}_0},
\end{equation}
where $(\cdot,\cdot)_{H^{-2},H^{2}_0}$ represents the usual duality product. 

\begin{definition}
Let $\eta$ be the weight function introduced at the beginning of Section \ref{weightedcontrol}.  For any $(y_0,y_1)\in X$ we define the minimal $L^2$-weighted control $v^\varepsilon(t)$ associated to (\ref{eq:sys_yeps}) as the function 
\begin{equation}
v^\varepsilon(t)=\eta (t)  \varphi^\varepsilon_{xx}(1,t)\in L^2(0,T)  \label{eq_osbeam}
\end{equation}
where $\varphi^\varepsilon$ is the solution of the adjoint system (\ref{eq:sys_pheps1d}) with initial data $(\varphi_0^\varepsilon,\varphi_1^{\varepsilon})$, the minimizer of 
\begin{equation} \label{eq_J}
J^\varepsilon (\varphi_0,\varphi_1)=\frac{\varepsilon}{2} \int_{\Sigma_0} \eta (t) |\varphi^\varepsilon_{xx}(1,t)|^2 \; dt -<(y_0,y_1),({\varphi_0},{\varphi_1})>_{X,X^*} , 
\end{equation}
in $(\varphi_0,\varphi_1)\in X^*$.
\end{definition}

The existence and uniqueness of the minimal $L^2$-weighted control is easily obtained when $T>2$ from the results in \cite{lions2}. In fact, the main ingredient is the observability inequality (\ref{eq:obsineqeps}) that provides the coercivity of $J^\varepsilon$ when the support of $\eta$ is an interval with length greater than 2.

We show below that, if the initial data $(y_0,y_1)$ are smooth, then the same is true for $(\varphi_0^{\varepsilon},\varphi_1^{\varepsilon})$, the minimizer of $J^\varepsilon$. In order to state the controllability result we introduce the scale Hilbert spaces associated to the associated operator. Let $A^\varepsilon_0:D(A^\varepsilon_0)\subset L^2(0,1)\to L^2(0,1)$ be the unbounded operator defined by $A^\varepsilon_0 = -\partial^2_{xx} +\varepsilon \partial^4_{xxxx}$ with domain $D(A^\varepsilon_0)=H^4\cap H^2_0(0,1)$. It is easy to see that $A^{\varepsilon}_0$ is a dissipative selftadjoint operator. 

We also define the unbounded skew-adjoint operator on $X=L^2(0,1)\times H^{-2}(0,1)$, 
$$
A^\varepsilon=\left( 
\begin{array}{cc}
0 & I \\
-A_0^\varepsilon & 0
\end{array}
\right), \qquad D(A^\varepsilon)=H^2_0(0,1)\times L^2(0,1). 
$$
Associated to $A^\varepsilon$ we consider the usual scale of Hilbert spaces $X_\alpha=D((A^\varepsilon)^\alpha)$, $\alpha >0$. Note that if we use the duality product (\ref{eq_dua_pr}) then 
$$
(A^\varepsilon)^*:D((A^\varepsilon)^*)\subset X^* \to X^*, 
$$
is given by 
$$
(A^\varepsilon)^*=\left( 
\begin{array}{cc}
0 & -I \\
A_0^\varepsilon & 0
\end{array}
\right), \qquad D((A^\varepsilon)^*)=X^*_1=X.
$$
In general, $X_\alpha^*=D(((A^\varepsilon)^*)^\alpha)=D((A^\varepsilon)^{\alpha+1})$.

The following result is a direct consequence of the results in \cite{ErvedozaZuazua}:

\begin{theorem} \label{th_EZ} 
Given any $(y_0,y_1)\in X=L^2(0,1)\times H^{-2}(0,1)$, there exists a unique weighted control $v^\varepsilon$ of system (\ref{eq:sys_yeps}) satisfying \eqref{eq_osbeam}.
This control is the one that minimizes the norm
$$
\int_0^T \frac{|v^\varepsilon|^2}{\eta(t)} dt.
$$
Furthermore, if $(y_0,y_1) \in D((A^\varepsilon)^\alpha)$ for some $\alpha>0$, the control $v^\varepsilon$ satisfies 
$$
v^\varepsilon \in H^{\alpha}_0(0,T) \bigcap_{k=0}^{[\alpha]} C^k([0,T]),
$$
and the corresponding $(\psi_0^{T,\varepsilon},\psi_1^{T,\varepsilon})\in X_\alpha^* = X_{\alpha+1}$. In particular, the controlled solution $y$ belongs to
$$
(y,y')\in C^\alpha([0,T];X_0)\bigcap_{k=0}^{[\alpha]} C^k([0,T];X_{\alpha-k}).
$$ 
Finally, the following estimate holds,
\begin{equation} \label{eq_est_cc}
\| v^\varepsilon\|_{H^\alpha_0(0,T)} \leq C \| (y_0,y_1)\|_{X_{\alpha}}.
\end{equation}
\end{theorem}
{\sc Proof.} The proof is a straightforward adaptation of one of the examples in \cite{ErvedozaZuazua} corresponding to the boundary controllability of the wave equation. The only difference is that here the operator $A_0^\varepsilon$ is a fourth order one with  different scale Hilbert spaces. \bfin

We finally observe that, for smooth solutions and controls the minimizer of (\ref{eq_J}),  $(\varphi_0^\varepsilon,\varphi_1^{\varepsilon})$ can be characterized as the solution of the following optimality system 
\begin{equation}
\label{eq_optsyseps1}
\left\{
\begin{aligned}
&\varphi_{tt}^\varepsilon + \varepsilon  \varphi^\varepsilon_{xxxx} -  \varphi^\varepsilon_{xx} =0, & \text{ in } Q_T, \\
&\varphi^\varepsilon (0,\cdot)=\varphi^\varepsilon_x (0,\cdot)=0, & \text{ in } (0,T)\\
&\varphi^\varepsilon (1,\cdot)=\varphi^\varepsilon_x (1,\cdot)=0, & \text{ in } (0,T)\\
&\varphi^\varepsilon(\cdot,0)=\varphi_0^{\varepsilon}, \quad \varphi^\varepsilon_t(\cdot,0)=\varphi_1^{\varepsilon}, & \text{ in } \Omega,
\end{aligned}
\right. 
\end{equation}
\begin{equation}
 \label{eq_optsyseps2}
 \left\{
\begin{aligned}
&y_{tt}^\varepsilon + \varepsilon  y^\varepsilon_{xxxx} -  y^\varepsilon_{xx} =0, & \text{ in } Q_T, \\
&y^\varepsilon (0,\cdot)=y^\varepsilon (1,\cdot)=y^\varepsilon_x (0,\cdot)=0, & \text{ in } (0,T)\\
&y^\varepsilon_x (1,\cdot)=\eta(\cdot)\varphi_{xx}^\varepsilon(1,\cdot), & \text{ in } (0,T)\\
&y^\varepsilon(\cdot,0)=y_0, \quad y^\varepsilon_t(\cdot,0)=y_1, & \text{ in } \Omega, \\
&y^\varepsilon(\cdot,T)= y^\varepsilon_t(\cdot,T)=0, & \text{ in }\Omega.
\end{aligned}
\right.
\end{equation}

\subsection{Minimal $L^2$-weighted controls for the wave equation}
 
Consider now the wave equation (\ref{eq:sys_y}) with initial data $(y_0,y_1)\in L^{2}\times H^{-1}$ and control $v\in L^2(0,T)$ for which 
\begin{equation} \label{eq yT0_w}
y (\cdot,T)=y_t(\cdot,T)=0 \quad \text{ in } (0,1).
\end{equation}

Let us introduce the homogeneous adjoint system to (\ref{eq:sys_y}):
\begin{equation} 
\label{eq_2we}
\left\{ 
\begin{aligned}
&\varphi_{tt} - \varphi_{xx} =0, & \text{ in } Q_T, \\
&\varphi (0,\cdot)=\varphi (1,\cdot)=0, & \text{ in } (0,T),\\
&\varphi (\cdot,T)=\varphi_0^T, \quad \varphi_t(\cdot,T)=\varphi_1^T, & \text{ in } \Omega.
\end{aligned}
\right.
\end{equation}

Let $Z=L^2(0,1)\times H^{-1}(0,1)$ and $Z^*=H^1_0(0,1)\times L^2(0,1)$ its dual, with duality product given by 
$$
<(y_0,y_1),(\varphi_0,\varphi_1)>_{Z,Z^*} = \int_\Omega y_0\varphi_1 \; dx - (y_1,\varphi_0)_{H^{-1},H^{1}_0},
$$
where $(\cdot,\cdot)_{H^{-1},H^{1}_0}$ represents the usual duality product. 

\begin{definition}
Let $\eta$ be the weight function introduced at the beginning of Section \ref{weightedcontrol}.  For any $(y_0,y_1)\in X$ we define the minimal $L^2$-weighted control $v(t)$ associated to (\ref{eq:sys_y}) as the function 
\begin{equation}
v(t)=\eta (t)  \varphi_{x}(1,t) \in L^2(0,T),\label{eq_oswave}
\end{equation} 
where $\varphi$ is the solution of the adjoint system (\ref{eq_2we}) with initial data the minimizer of the functional 
\begin{equation} \label{eq_Jwe}
J (\varphi_0,\varphi_1)=\frac{1}{2} \int_{0}^T \eta (t) | \varphi_x(1,t)|^2 \; dt -<(y_0,y_1),({\varphi}(x,0),{\varphi}_t(x,0))>_{Z,Z^*} , 
\end{equation}
in $(\varphi_0,\varphi_1)\in Z^*$.
\end{definition}

The existence and uniqueness of the minimal $L^2$-weighted control is well-known when $T>2$. The main ingredient is the following well-known observability inequality that provides the coercivity of $J$ when the support of $\eta$ is an interval with length greater than 2,
\begin{equation} \nonumber
\| \varphi_0 \|_{H^1_0(0,1)}^2+ \| \varphi_1 \|_{L^2(0,1)}^2  \leq C\int_0^T \left| \varphi_{x}(1,t)\right|^2 dt, 
\end{equation}
when $T>2$, for some constant $C>0$.

As for the beam model, regularity of the controls can be improved if we have smoother initial conditions. We define the unbounded skew-adjoint operator on $Z=L^2(0,1)\times H^{-1}(0,1)$, 
$$
A=\left( 
\begin{array}{cc}
0 & I \\
-\partial^2_{xx} & 0
\end{array}
\right), \qquad D(A)=H^1_0(0,1)\times L^2(0,1). 
$$
Associated to $A$ we consider the usual scale of Hilbert spaces $Z_\alpha=D(A^\alpha)$, $\alpha >0$. Note that in this case $Z_\alpha^*=D((A^*)^\alpha)=D(A^{\alpha+1})$. The following result is proved in \cite{ErvedozaZuazua},

\begin{theorem} \label{th_EZwe}
Given any $(y_0,y_1)\in Z=L^2(0,1)\times H^{-1}(0,1)$, there exists a unique weighted control $v$ of system (\ref{eq:sys_y}) satisfying \eqref{eq_oswave}).
This control is the one that minimizes the norm
$$
\int_0^T \frac{|v(t)|^2}{\eta(t)} dt.
$$
Furthermore, if $(y_0,y_1) \in D(A^\alpha)$ for some $\alpha>0$, the control $v$ satisfies 
$$
v \in H^{\alpha}_0(0,T) \bigcap_{k=0}^{[\alpha]} C^k([0,T])
$$
and the corresponding $(\psi_0^{T,\varepsilon},\psi_1^{T,\varepsilon})\in Z_\alpha^* = Z_{\alpha+1}$. In particular, the controlled solution $y$ belongs to
$$
(y,y')\in C^\alpha([0,T];X_0)\bigcap_{k=0}^{[\alpha]} C^k([0,T];X_{\alpha-k}).
$$ 
\end{theorem}

For smooth solutions the minimizer of (\ref{eq_Jwe}),  $(\varphi_0,\varphi_1)$ is characterized as the solution of the following optimality system 
\begin{equation}
 \label{eq_optsyswe1}
 \left\{
\begin{aligned}
&\varphi_{tt}  -  \varphi_{xx} =0, & \text{ in } Q_T, \\
&\varphi (0,\cdot)=\varphi_x (0,\cdot)=0, & \text{ in } (0,T),\\
&\varphi(\cdot,0)=\varphi_0^{T}, \quad \varphi_t(\cdot,0)=\varphi_1^{T}, & \text{ in } \Omega,
\end{aligned}
\right. 
\end{equation}
\begin{equation}
\label{eq_optsyswe2}
 \left\{
\begin{aligned}
& y_{tt}  -  y_{xx} =0, & \text{ in } Q_T, \\
& y (0,\cdot)=0, \quad  y (1,\cdot)=\eta(\cdot)\varphi_{x}(1,\cdot),& \text{ in } (0,T),\\
& y(\cdot,0)=y_0, \quad y_t(\cdot,0)=y_1, & \text{ in } \Omega, \\
& y(\cdot,T)= y_t(\cdot,T)=0, & \text{ in } \Omega.
\end{aligned}
\right.
\end{equation}

\section{Asymptotic expansion of the direct and adjoint systems}\label{asymptotic_expansion}

We construct in this section an asymptotic expansion for the solution of the following general system 
\begin{equation} 
\label{eq_ep_gen}
\left\{
\begin{aligned}
& y_{tt}^\varepsilon + \varepsilon  y^\varepsilon_{xxxx} - y^\varepsilon_{xx} =0, & \mbox{in}\quad  Q_T \\
& y^\varepsilon (0,\cdot)=y^\varepsilon (1,\cdot)=y^\varepsilon_{x} (0,\cdot)=0, \quad y^\varepsilon_{x} (1,\cdot)=v^\varepsilon, & \text{ in }\quad  (0,T),\\
& (y^\varepsilon(\cdot,0),y^\varepsilon_t(\cdot,0))=(y_0^{\varepsilon},y_1^\varepsilon), & \mbox{in}\quad  \Omega.
\end{aligned}
\right.
\end{equation}
This will be used for both the direct problem (\ref{eq:sys_yeps}), for which the initial data will not depend on $\varepsilon$, and the adjoint one for which $v^\varepsilon=0$.   

Of course, this requires some a priori information on the asymptotics for the control $v^\varepsilon$ and the initial data. In view of Theorem \ref{theoremLions}, we assume that the function $v^{\eps}$ is in the form 
\begin{equation} \label{exp_cont}
\sqrt{\eps} v^\varepsilon=\sum_{k=0}^N \varepsilon^{k/2} v^k, 
\end{equation}
the functions $v^k$, $k\geq 0$ being known. Note that, as we are dealing with controls $v^\varepsilon$ which vanish in a neighborhood of $t=0$, we assume the same for $v^k$, i.e. 
\begin{equation} \label{cond_contr}
v^k(t)=0, \quad \mbox{ in a neighborhood of $t=0$, for all $k\geq 0$.} 
\end{equation}

Concerning the initial data, we do not make assumptions for the moment. Solutions of (\ref{eq_ep_gen}) have two boundary layers at $x=0,1$ respectively. In order to have a convergent result for the asymptotics we will require a similar behavior for the initial data that we make precise later. 

\subsection{Formal asymptotic expansion}

In order to construct an asymptotic expansion of $y^{\eps}$, we use the method of matched asymptotic expansion (see \cite{ColeBook96,eckhausbook79}). Let us consider the following formal outer and inner expansions
\begin{eqnarray}
&& y^\varepsilon (x,t)=y^\varepsilon_{out} (x,t) \sim \sum_{k=0}^N \varepsilon ^{k/2} y^k (x,t), \quad x\in(0,1), \; t\in (0,T),  \label{eq outer}\\
&& y^\varepsilon (x,t)=y^\varepsilon_{in} (z,t) \sim \sum_{k=0}^N \varepsilon ^{k/2} Y^k (z,t), \quad z=\frac{x}{\sqrt{\varepsilon}}\in(0,\varepsilon^{-1/2}), \; t\in (0,T), \label{eq inner1}\\
&& y^\varepsilon (x,t)=y^\varepsilon_{in} (w,t) \sim \sum_{k=0}^N \varepsilon ^{k/2} S^k (w,t), \quad w=\frac{1-x}{\sqrt{\varepsilon}}\in(0,\varepsilon^{-1/2}), \; t\in (0,T). \label{eq inner2} 
\end{eqnarray}
Here the boundary layer (inner region) occurs near $x=0$ and $x=1$ and it is of ${\cal O} (\sqrt{\varepsilon})$ size, and the outer region is the subset of $(0,1)$ consisting of the points far from the boundary layer, it is of  ${\cal O} (1)$ size. In (\ref{eq outer})-(\ref{eq inner2}) we have noted by $y^\eps_{out}$ and $y^\eps_{in}$ these outer and inner expansion respectively.
There is an intermediate region between them with size ${\cal O} (\varepsilon^\gamma)$, $\gamma\in(0,1)$ to be determined. To construct the approximate solution we require that the inner expansion close to $x=0$, given by (\ref{eq inner1}), equals the outer expansion (\ref{eq outer}) in the intermediate region $\varepsilon^{1/2} << x<< 1$. Analogously, the inner expansion close to $x=1$, given by (\ref{eq inner2}), must coincide with the outer expansion (\ref{eq outer}) in the intermediate region $0 << x<< 1-\varepsilon^{1/2}$.  These conditions are the so-called matching asymptotic conditions. The solution in this intermediate region will be noted as $y_{match}$ below. 

Inserting (\ref{eq outer}) into equation (\ref{eq_ep_gen}) and making equal the terms with the same power in $\varepsilon$ yields,
\begin{equation} \label{eq_asx}
\begin{array}{rcl}
\varepsilon^0&:& y^0_{tt}-y^0_{xx}=0, \\
\varepsilon^{1/2}&:& y^1_{tt}-y^1_{xx}=0, \\
\varepsilon^{1}&:& y^2_{tt}-y^2_{xx}=-y^0_{xxxx}, \\
\varepsilon^{3/2}&:& y^3_{tt}-y^3_{xx}=-y^1_{xxxx}, \\
...&:& ...
\end{array}
\end{equation}
 Similarly, for the inner expansions we obtain 
\begin{equation} \label{eq_asY}
\begin{array}{rcl}
\varepsilon^{-1}&:& Y^0_{zzzz}-Y^0_{zz}=0, \\
\varepsilon^{-1/2}&:& Y^1_{zzzz}-Y^1_{zz}=0, \\
\varepsilon^{0}&:& Y^2_{zzzz}-Y^2_{zz}=-Y^0_{tt}, \\
\varepsilon^{1/2}&:& Y^3_{zzzz}-Y^3_{zz}=-Y^1_{tt}, \\
...&:& ...
\end{array}
\end{equation}
and
\begin{equation} \label{eq_asS}
\begin{array}{rcl}
\varepsilon^{-1}&:& S^0_{wwww}-S^0_{ww}=0, \\
\varepsilon^{-1/2}&:& S^1_{wwww}-S^1_{ww}=0, \\
\varepsilon^{0}&:& S^2_{wwww}-S^2_{ww}=-S^0_{tt}, \\
\varepsilon^{1/2}&:& S^3_{wwww}-S^3_{ww}=-S^1_{tt}, \\
...&:& ...
\end{array}
\end{equation}
We impose the boundary conditions to the inner solution. This reads,
\begin{equation} \label{eq_bo_as}
Y^k(0,t)=Y^k_z(0,t)=S^k(0,t)=0,\qquad S^k_w(0,t)=-v^k(t), \quad t\in(0,T).
\end{equation}
Therefore the solutions $Y_k$ and $S_k$ can be computed up to some integration constants. For $k=0,1$ we easily obtain,
\begin{eqnarray*}
Y^0(z,t) &=& C_{0,1}(t) (e^{-z}+z-1)+C_{0,2}(t) (e^{z}-z-1), \\
S^0(w,t) &=& D_{0,1}(t) (e^{-w}+w-1)+D_{0,2}(t) (e^{w}-w-1)-wv^0(t), \\
Y^1(z,t) &=& C_{1,1}(t) (e^{-z}+z-1)+C_{1,2}(t) (e^{z}-z-1), \\
S^1(w,t) &=& D_{1,1}(t) (e^{-w}+w-1)+D_{1,2}(t) (e^{w}-w-1)-wv^1(t).   
\end{eqnarray*}
We determine these constants by asymptotically matching the inner and outer solutions. 
The match consists of requiring that the intermediate limits ($\varepsilon \to 0^+, \; x\to 0^+, z=x/\varepsilon^{1/2} \to \infty$, $w=(1-x)/\varepsilon^{1/2} \to \infty$) of the inner and outer solutions agree. 
We first remark that 
$$
C_{k,2}(t)= D_{k,2}(t)=0, \quad t\geq 0, \quad k=0,1,
$$ 
since otherwise it is not possible to match the exponential growing behavior for $z,w\to \infty$.

The rest of the matching conditions are obtained by considering the Taylor expansion of the outer solution $y^\varepsilon_{out}(x,t)$ at the inner region. Thus, for $x<<1$ we have formally
\begin{eqnarray} \nonumber
y^\varepsilon_{out}(x,t)&=&y^0(0,t) +xy^0_x(0,t)+\frac{x^2}{2} y^0_{xx}(0,t) + \ldots \\
\nonumber
&&+ \varepsilon^{1/2}\left( y^1(0,t) +xy^1_x(0,t)+\frac{x^2}{2} y^1_{xx}(0,t) + \ldots\right) \\
\nonumber
&&+ \varepsilon\left( y^2(0,t) +xy^2_x(0,t)+\frac{x^2}{2} y^2_{xx}(0,t) + \ldots\right)+ \ldots\\
\nonumber
&=&y^0(0,t) +\varepsilon^{1/2} \left( zy^0_x(0,t) +  y^1(0,t)\right) +\varepsilon \left( \frac{z^2}{2} y^0_{xx}(0,t) + z y^1_x(0,t) + y^2(0,t) \right) \\
\label{eq_as_ex}
&&+ \varepsilon^{3/2} \left( \frac{z^3}{6} y^0_{xxx}(0,t) + \frac{z^2}{2} y^1_{xx}(0,t) + z y^2_x(0,t)+ y^3(0,t)\right) + \ldots
\end{eqnarray}
that we match with the inner expansion (\ref{eq inner1}), i.e. $y_{out}^\varepsilon (x,t)=y^\varepsilon_{in}(z,t)$, as $z=x/\varepsilon^{1/2} \to \infty$. 

Analogously for $1-x<<1$ we have
\begin{eqnarray} \nonumber
y^\varepsilon_{out}(x,t)&=&y^0(1,t) +\varepsilon^{1/2} \left( -wy^0_x(1,t) +  y^1(1,t)\right) \\
\nonumber
&&+\varepsilon \left( \frac{w^2}{2} y^0_{xx}(1,t) - w y^1_x(1,t) + y^2(1,t) \right) \\
\label{eq_as_ex_w}
&&+ \varepsilon^{3/2} \left( -\frac{w^3}{6} y^0_{xxx}(1,t) + \frac{w^2}{2} y^1_{xx}(1,t) - w y^2_x(1,t)+ y^3(1,t)\right) + \ldots
\end{eqnarray}
that we match with the inner expansion (\ref{eq inner2}), i.e. $y_{out}^\varepsilon (x,t)=y^\varepsilon_{in}(w,t)$, as $w=(1-x)/\varepsilon^{1/2} \to \infty$. 

\par\noindent
Thus, equaling terms with the same power in $\varepsilon$ we obtain the different matching conditions: 

\par\noindent
$\bullet$ Order $\varepsilon^0$-  The leading order term is given by 
\begin{equation} \label{eq_mat01}
y^0(0,t)=\lim_{z\to \infty} Y^0(z,t)=\lim_{z\to \infty} C_{0,1}(t) \left( e^{-z}+z-1 \right).
\end{equation} 
This implies $C_{0,1}(t)=0$ and therefore,
\begin{equation} \label{eq_mat1}
y^0(0,t)=0, \quad  Y^0(z,t)=0, \quad t\geq 0.
\end{equation}
In this way, we obtain a boundary condition for the zero order term of the outer expansion $y^\eps_{out}$ at $x=0$ and the zero order term for the inner expansion $y^\eps_{in}$ in (\ref{eq inner1}).  Observe that the match occurs in the intermediate region where $1<<z=x/\varepsilon^{1/2}$ as well as $x<<1$. Thus, the size of the overlap region is  $\varepsilon^{1/2}<<x<<1$. 

An analogous argument for the solution near $x=1$ provides 
\begin{equation} \label{eq_mat02}
y^0(1,t)=\lim_{w\to \infty} S^0(w,t)=\lim_{w\to \infty} \left(D_{0,1}(t) \left( e^{-w}+w-1 \right)-w\,v^0(t) \right).
\end{equation}
This implies $D_{0,1}(t)=v^0(t)$ and therefore,
\begin{equation} \label{eq_mat2}
y^0(1,t)=-v^0(t), \quad S^0(w,t)=-y^0(1,t)\left( e^{-w}-1 \right), \quad t\geq 0.
\end{equation}
Again, this provides a boundary condition for the zero order term of the outer expansion $y_{out}$ at $x=1$ and the zero order term for the inner expansion $y_{in}$ in (\ref{eq inner2}).

The first equation in (\ref{eq_asx}), together with the boundary conditions in (\ref{eq_mat1}) and (\ref{eq_mat2}), determine $y^0$ up to some initial conditions that we write $(y_0^{0},y_1^{0})$. In this way, $y^0$ is the solution of the following system 
\begin{equation} \label{eq_sis0}
\left\{ \begin{array}{ll}
y^0_{tt} - y^0_{xx} =0, & \mbox{in}\, Q_T,  \\
y^0 (0,\cdot)=0, \quad y^0 (1,\cdot)= -v^0, & \mbox{in}\, (0,T), \\
y^0(\cdot,0)=y_0^{0}, \quad y_{t}^0(\cdot,0)=y_1^{0}, &  \mbox{in}\, (0,1).
\end{array} \right.
\end{equation}
On the other hand, $Y^0$ and $S^0$ are given in (\ref{eq_mat1}) and (\ref{eq_mat2}) respectively. Therefore we have computed the zero order terms of both $y_{out}$ and $y_{in}$ in (\ref{eq outer})-(\ref{eq inner1}).

As usual, a so-called composite approximation over the complete interval $x\in(0,1)$ is obtained by adding, at each order, the inner and outer expansion and then by subtracting their common part,
\begin{eqnarray} \nonumber
y^\varepsilon(x,t)&=&y_{out}^\varepsilon (x,t)+y^\varepsilon_{in}(z,t)+y^\varepsilon_{in}(w,t)- y_{match}^\varepsilon (x,t) +{\cal O}(\varepsilon^{1/2})\\ \label{sol_ord0}
&=& y^0(x,t)+v^0(t) e^{-w}+{\cal O}(\varepsilon^{1/2})=y^0(x,t)-y^0(1,t) e^{-w}+{\cal O}(\varepsilon^{1/2}),  
\end{eqnarray}
where the common part, named here $y_{match}^\varepsilon (x,t)$, is either the outer or inner solution in the matching region. In this case, $y_{match}^\varepsilon (x,t)=0$ in the neighborhood of $x=0$ and $y_{match}^\varepsilon (x,t)=-v^0(t)$ in the neighborhood of $x=1$.

Finally, we state the hypotheses on the initial data for $y^\varepsilon$. In order to have the asymptotics in (\ref{sol_ord0}), we should have  
\begin{eqnarray*}
y^\varepsilon(x,0)&=&y^0(x,0)-y^0(1,0) e^{-(1-x)/\varepsilon^{1/2}}+{\cal O}(\varepsilon^{1/2})= y^0_0(x)+{\cal O}(\varepsilon^{1/2}),\\
y^\varepsilon_t(x,0)&=&y^0_t(x,0)-y^0_t(1,0) e^{-(1-x)/\varepsilon^{1/2}}+{\cal O}(\varepsilon^{1/2})= y^0_1(x)+{\cal O}(\varepsilon^{1/2}),
\end{eqnarray*}
where we have used that $y^0(1,0)=v^0(0)=0$ and $y^0_t(1,0)=v^0_t(0)=0$ in view of (\ref{cond_contr}). According to this, our assumption on the initial data is that it can be written as follows:
$$
y^\varepsilon(x,0)=y^0_0(x)+{\cal O}(\varepsilon^{1/2}), \quad y^\varepsilon_t(x,0)=y^0_1(x)+{\cal O}(\varepsilon^{1/2}),
$$
for some known functions $(y^0_0,y^0_1)$. 

\par\noindent
$\bullet$ Order $\varepsilon^{1/2}$- Next we match to first order in (\ref{eq_as_ex}), keeping terms of order $\varepsilon^{1/2}$, with the inner solution $y^\varepsilon_{in}$ in (\ref{eq inner1}). Then, we have
\begin{equation} \label{eq_mat01_b}
\lim_{z\to \infty} \left( zy^0_x(0,t)+y^1(0,t) \right) =\lim_{z\to \infty} Y^1(z,t)=\lim_{z\to \infty} C_{1,1}(t) \left( e^{-z}+z-1 \right).
\end{equation} 
This implies $C_{1,1}(t)=y^0_x(0,t)$ and therefore,
\begin{equation} \label{eq_mat1_b}
y^1(0,t)=-y^0_x(0,t), \quad Y^1(z,t)=-y^1(0,t)\left( e^{-z}+z-1 \right), \quad t\geq 0.
\end{equation}
Since matching now requires that $x^2<<\varepsilon^{1/2}$ the size of the overlap region is smaller than it was for the leading order match. Its extents is $\varepsilon^{1/2}<<x<<\varepsilon^{1/4}$ ($\varepsilon \to 0^+$). 
 
Proceeding in a similar way for the inner solution near $x=1$, the order $\varepsilon^{1/2}$ for $y^\eps_{out}$ in (\ref{eq_as_ex_w}) must be equal to the corresponding for $y^\eps_{in}$ in (\ref{eq inner2}). We easily obtain 
\begin{equation} \label{eq_mat01_b2}
\lim_{w\to \infty} \left( -wy^0_x(1,t)+y^1(1,t) \right) =\lim_{w\to \infty} S^1(w,t)=\lim_{w\to \infty} \left( D_{1,1}(t) \left( e^{-w}+w-1 \right) - wv^1(t) \right).
\end{equation} 
This implies $D_{1,1}(t)=-y^0_x(1,t)+v^1(t)$ and therefore,
\begin{equation} 
\label{eq_mat1_b2}
y^1(1,t)=y^0_x(1,t)-v^1(t), \quad S^1(w,t)=-y^1(1,t)\left( e^{-w}+w-1 \right)- wv^1(t), \quad t\geq 0.
\end{equation}
In particular $y^1$ is determined, up to some initial data that we write $(y_0^1,y_1^1)$,  as the solution of the following system,
\begin{equation} \label{eq_sis11}
\left\{ \begin{array}{ll}
y^{1}_{tt} - y^{1}_{xx} =0, & \mbox{in}\, Q_T  \\
y^1 (0,\cdot)=- y_x^0(0,\cdot), \quad y^1 (1,\cdot)= y_x^0(1,\cdot)-v^1, & \mbox{in}\, (0,T), \\
y^1(\cdot,0)=y_0^{1}, \quad y_{t}^1(\cdot,0)=y_1^{1}, &  \mbox{in}\, (0,1).
\end{array} \right.
\end{equation}
Once the boundary layer solution is determined, we recover  the following approximation in the whole interval $x\in(0,1)$,
\begin{eqnarray} \nonumber
&&y^\varepsilon(x,t)=y_{out}^\varepsilon (x,t)+y^\varepsilon_{in}(z,t)+y^\varepsilon_{in}(w,t)- y_{match}^\varepsilon (x,t) +{\cal O}(\varepsilon)\\ \nonumber
&&\quad =\sum_{k=0}^1 \varepsilon^{k/2}\left( y^k(x,t)+Y^k(z,t)+S^k(w,t)\right) - y_{match}^\varepsilon (x,t) +{\cal O}(\varepsilon)\\
 \label{eq_sec_or}
&&\quad = \sum_{k=0}^1 \varepsilon^{k/2} \left[ y^k(x,t)-
y^{k}(0,t) e^{-z} - y^{k}(1,t) e^{-w} \right] +{\cal O}(\varepsilon),
\end{eqnarray}
with 
$$
y_{match}^\varepsilon (x,t)=-v^0(t)+ \varepsilon^{1/2} \left( y^0_x(0,t)(z-1)-y^0_x(1,t)(w-1)-v^1(t)\right).
$$ 

We now state the hypotheses on the asymptotic expansion of the initial data in order to have (\ref{eq_sec_or}). Note that, 
\begin{eqnarray*}
y^\varepsilon(x,0)&=&\sum_{k=0}^1 \varepsilon^{k/2} \left[ y^k(x,0)-
y^{k}(0,0) e^{-z} - y^{k}(1,0) e^{-w} \right] +{\cal O}(\varepsilon)\\
&=& \sum_{k=0}^1 \varepsilon^{k/2} \left[ y^k_0(x)-
y^{k}_0(0) e^{-z} - y^{k}_1(0) e^{-w} \right] +{\cal O}(\varepsilon),
\end{eqnarray*}
and an analogous formula holds for $y^\varepsilon_t(x,0)$.
Therefore, we assume
\begin{eqnarray*}
y^\varepsilon(x,0)&=&\sum_{k=0}^1 \varepsilon^{k/2} \left[ y^k_0(x)-
y^{k}_0(0) e^{-x/\varepsilon^{1/2}} - y^{k}_0(1) e^{-(1-x)/\varepsilon^{1/2}} \right] +{\cal O}(\varepsilon),\\
y^\varepsilon_t(x,0)&=&\sum_{k=0}^1 \varepsilon^{k/2} \left[ y^k_1(x)-
y^{k}_1(0) e^{-x/\varepsilon^{1/2}} - y^{k}_1(1) e^{-(1-x)/\varepsilon^{1/2}} \right] +{\cal O}(\varepsilon),
\end{eqnarray*}
for some $(y_0^k,y_1^k)$ with $k=0,1$.
Note that in these expressions some terms vanish, as $y^0_0(0), \; y^0_1(0)$, but others not necessarily, as $y^1_0(0), \; y_1^1(0)$. However, we keep this notation to have more abriged formulas. 

\par\noindent
$\bullet$ Order $\varepsilon^{1}$- Next we match to second order in (\ref{eq_as_ex}) with $y^\varepsilon_{in}$, keeping terms of order $\varepsilon$. Then, we have
\begin{equation} \label{eq_mat01_b3}
\lim_{z\to \infty} \left( \frac{z^2}{2}y^0_{xx}(0,t)+zy^1_x(0,t)+y^2(0,t) \right) =\lim_{z\to \infty} Y^2(z,t)=\lim_{z\to \infty} C_{2,1}(t) \left( e^{-z}+z-1 \right),
\end{equation} 
where we have solved the equation for $Y^2$ in (\ref{eq_asY}) taking into account that $Y^0=0$ and the boundary conditions in (\ref{eq_bo_as}). This implies 
$$
y^0_{xx}(0,t)=0, \quad C_{2,1}(t)=y^1_x(0,t).
$$ 
The first condition is satisfied as long as $y^0(x,t)$ is a sufficiently smooth solution of (\ref{eq_sis0}) to have this trace. Concerning the second condition, this implies
\begin{equation} \label{eq_mat1_c3}
y^2(0,t)=-y^1_x(0,t), \quad Y^2(z,t)=-y^2(0,t)\left( e^{-z}+z-1 \right), \quad t\geq 0.
\end{equation}
Since matching now requires that $x^3<<\varepsilon$ the size of the overlap region is smaller than it was for the first order match. Its extents is $\varepsilon^{1/2}<<x<<\varepsilon^{1/3}$ ($\varepsilon \to 0^+$). 

Proceeding in a similar way for the inner solution near $x=1$ we easily obtain 
\begin{equation} \label{eq_mat01_c3}
\lim_{w\to \infty} \left( \frac{w^2}{2}y^0_{xx}(1,t)-wy^1_x(1,t)+y^2(1,t) \right) =\lim_{w\to \infty} S^2(w,t),
\end{equation} 
where $S^2(w,t)$ is solution of 
$$
\left\{
\begin{array}{ll}
S_{wwww}^2-S_{ww}^2=-S^0_{tt}=-v^0_{tt}(t)(e^{-w}-1), & w\in (0,\varepsilon^{-1/2}), \; t>0,\\
S^2(0,t)=0, \quad S^2_w(0,t)=-v^2(t), & t>0.
\end{array} 
\right.
$$
The general solution is given by
$$
S^2(w,t)=v^0_{tt}(t)\frac{w}{2} (e^{-w}-w-1)+D_{2,1}(t)(e^{-w}+w-1)+D_{2,2}(t)(e^{w}-w-1)-wv^2(t). 
$$
Obviously, $D_{2,2}(t)$ must be zero to satisfy the matching condition (\ref{eq_mat01_c3}). Moreover, we must also have 
$$
y^0_{xx}(1,t)=-v^0_{tt}(t), \qquad D_{2,1}(t) = - y^1_x(1,t)+v^2(t)+\frac12 v^0_{tt}(t).
$$
The first condition is satisfied as long as $y^0$ is a sufficiently smooth solution of (\ref{eq_sis0}). The second identity gives
\begin{equation} \label{eq_mat1_c}
\begin{aligned}
&y^2(1,t)= y^1_x(1,t)-v^2(t)-\frac12 v^0_{tt}(t), \\ 
&S^2(w,t)=v^0_{tt}(t)\frac{w}{2} (e^{-w}-w-1)-y^2(1,t)(e^{-w}+w-1) - w\,v^2(t).
\end{aligned}
\end{equation}
In particular, $y^2$ is determined, up to some initial data $(y_0^{2}, y_1^{2})$, by the solution of the following system,
\begin{equation} \label{eq_sis11_2}
\left\{ \begin{array}{ll}
y^{2}_{tt} - y^{2}_{xx} =-y^0_{xxxx}, & \mbox{in}\,\, Q_T, \\
y^2 (0,\cdot)=- y_x^1(0,\cdot), \quad y^2 (1,\cdot)=  y^1_x(1,\cdot)-v^2-\frac12 v^0_{tt}, & \mbox{in}\,\, (0,T), \\
y^2(\cdot,0)=y_0^{2}, \quad y_{t}^2(\cdot,0)=y_1^{2}, &  \mbox{in}\,\, (0,1).
\end{array} \right.
\end{equation}

Once the boundary layer solution is determined, one may construct the uniform approximation in the whole interval $x\in(0,1)$,
\begin{equation} 
\begin{aligned}
y^\varepsilon(x,t)=&y_{out}^\varepsilon (x,t)+y^\varepsilon_{in}(z,t)+y^\varepsilon_{in}(w,t)- y_{match}^\varepsilon (x,t) +{\cal O}(\varepsilon^{3/2}) \\
 =&\sum_{j=0}^2 \varepsilon^{j/2}\left( y^j(x,t)+Y^j(z,t)+S^j(w,t)\right) - y_{match}^\varepsilon (x,t) +{\cal O}(\varepsilon^{3/2})\\
 =& \sum_{j=0}^2 \varepsilon^{j/2} \left[ y^j(x,t)-
y^{j}(0,t) e^{-z} -\left( y^{j}(1,t)+\frac{w}{2} y^{j-2}_{tt}(1,t)\right) e^{-w} \right]+{\cal O}(\varepsilon^{3/2}). \label{eq_sec_or_c}
\end{aligned}
\end{equation}
Here we have assumed $y^{-2}=y^{-1}=0$. The assumptions on the initial data are now 
\begin{eqnarray*}
y^\varepsilon(x,0)&=&\sum_{k=0}^2 \varepsilon^{k/2} \left[ y^k_0(x)-
y^{k}_0(0) e^{-x/\varepsilon^{1/2}} - y^{k}_0(1) e^{-(1-x)/\varepsilon^{1/2}} \right] +{\cal O}(\varepsilon),\\
y^\varepsilon_t(x,0)&=&\sum_{k=0}^2 \varepsilon^{k/2} \left[ y^k_1(x)-
y^{k}_1(0) e^{-x/\varepsilon^{1/2}} - y^{k}_1(1) e^{-(1-x)/\varepsilon^{1/2}} \right] +{\cal O}(\varepsilon),
\end{eqnarray*}
for some $(y_0^k,y_1^k)$ with $k=0,1,2$.

\subsection{Convergence results}

We now check that the composite expressions we have determined are indeed uniform approximations of the solution $y^{\eps}$ of \eqref{eq_ep_gen} for $\eps$ small. 
We shall then apply this result both to the solution of the controlled problem (\ref{eq:sys_yeps}) and the adjoint one (\ref{eq:sys_pheps1d}). Note that system (\ref{eq:sys_yeps}) corresponds to the particular case $(y^\varepsilon_0,y^\varepsilon_1)=(y_0,y_1)$, independent of $\varepsilon$, while the adjoint system has $v^\varepsilon=0$.

We state below the convergence results for the three first terms in the asymptotic expansion (\ref{expansion_control}). We observe that each new term requires more regularity assumptions on the asymptotics of the initial data and the control. The proofs are technical and left to the appendix at the end of this paper. 
\begin{proposition} \label{pr_1}
Assume that $(y_0^{\varepsilon}, y_1^{\varepsilon}) \in X_1=H^2_0\times L^2$ and  consider $(y_0^{0}, y_1^{0}) \in Z_4\subset H^4\times H^3$,  $ v^\varepsilon\in H^2_0(0,T) $ and $v^0\in H^3_0(0,T)$  that satisfies the following estimates:
\begin{equation} \label{eq_l31_1}
\begin{aligned}
& \| (y_0^{\varepsilon}, y_1^{\varepsilon})-(y_0^{0}, y_1^{0}) \|_{H^1(0,1)\times L^2(0,1)} +\varepsilon^{1/2} \| y_{0,xx}^{\varepsilon}-y_{0,xx}^{0}\|_{L^2(0,1)}\leq C\varepsilon^{1/2},\\ 
& \| \varepsilon^{1/2} v^\varepsilon -v^0 \|_{H^2_0(0,T)} \leq C\varepsilon^{1/2}.
\end{aligned}
\end{equation}
Let $y^\varepsilon$ be the solution of (\ref{eq_ep_gen}) and $y^0$ be the solution of (\ref{eq_sis0}). Then, there exists a constant $C>0$, independent of $\varepsilon$, such that if we define the zero order approximation as follows
\begin{equation} \label{eq_zero_om}
y^{\varepsilon,0}(x,t)=y^0(x,t) +v^0(t) e^{-(1-x)/\varepsilon^{1/2}},
\end{equation}
the following estimates hold,
\begin{equation} 
\label{eq_l31_4}
\begin{aligned}
& \| y^\varepsilon-y^{\varepsilon,0}\|_{L^\infty(0,T;H^1_0(0,1))} + \| y^\varepsilon_t-y^{\varepsilon,0}_t \|_{L^\infty(0,T;L^2(0,1))} + \varepsilon^{1/2} \| y^\varepsilon_{xx}- y^{\varepsilon,0}_{xx} \|_{L^\infty(0,T;L^2(0,1))}\leq C  \varepsilon^{1/2} ,\\ 
&\| \varepsilon^{1/2}y^\varepsilon_{xx}(0,\cdot)- y^0_x(0,\cdot) \|_{L^2(0,T)} \leq C\varepsilon^{1/2}, \\ 
& \| \varepsilon^{1/2}y^\varepsilon_{xx}(1,\cdot)+ y^0_x(1,\cdot) \|_{L^2(0,T)} \leq C\varepsilon^{1/2}. 
\end{aligned}
\end{equation}
\end{proposition}

\begin{proposition} \label{pr_2}
Assume that $(y_0^{\varepsilon}, y_1^{\varepsilon}) \in X_1=H^2_0\times L^2$ and  consider, for $j=0,1$, $(y_0^{j}, y_1^{j}) \in Z_{5-j}\subset H^{5-j}\times H^{4-j}$, $ v^\varepsilon \in H^2_0(0,T)$ and $v^j \in H^{5-j}_0(0,T)$  that satisfies the following compatibility conditions 
\begin{equation}
 \label{eq_compo1}
 \begin{aligned}
& -y^0_{0,x} (0)=y_{0}^{1} (0), \quad y_{0,x}^{0} (1)=y_{0}^{1} (1), \\ 
& -y_{0,xxx}^{0} (0)=y_{0,xx}^{1} (0), \quad y_{0,xxx}^{0} (1)=y_{0,xx}^{1} (1).
\end{aligned}
\end{equation}
Let $y^\varepsilon$ be the solution of (\ref{eq_ep_gen}) and $y^0, \; y^1$ be the solutions of (\ref{eq_sis0}) and (\ref{eq_sis11}) with initial data $(y_0^{j}, y_1^{j})$, $j=0,1$, respectively. Define, 
\begin{eqnarray} \nonumber 
&&y^{\varepsilon,1}(x,t)= \sum_{j=0}^1 \varepsilon^{j/2} \left[ y^j(x,t)-
y^{j}(0,t) e^{-x/\varepsilon^{1/2}} -y^{j}(1,t) e^{-(1-x)/\varepsilon^{1/2}} \right] .
\label{eq_app1}
\end{eqnarray}
and assume that 
\begin{equation}\label{eq_estdi1}
\begin{aligned}
& \| (y_0^{\varepsilon}, y_1^{\varepsilon})-(y_0^{\eps,1}, y_1^{\eps,1}) \|_{H^1\times L^2} +\varepsilon^{1/2} \| y_{0,xx}^{\varepsilon}-y_{0,xx}^{\eps,1}\|_{L^2(0,1)}\leq C\varepsilon,\\ 
&\left\| \varepsilon^{1/2}v^\varepsilon - \sum_{j=0}^1 \varepsilon^{j/2} v^j\right\|_{H^2(0,T)} \leq C\varepsilon.
\end{aligned}
\end{equation}

Then, there exists a constant $C>0$, independent of $\varepsilon$, such that 
\begin{equation}
 \label{est_po1}
 \begin{aligned}
& \| y^\varepsilon-y^{\varepsilon,1} \|_{L^\infty(0,T;H^1_0)} + \| y^\varepsilon_t-y^{\varepsilon,1}_t \|_{L^\infty(0,T;L^2)} + \varepsilon^{1/2} \| y^\varepsilon_{xx}- y^{\varepsilon,1}_{xx} \|_{L^\infty(0,T;L^2)}\leq C  \varepsilon^{3/4},\\ 
& 
\| \varepsilon^{1/2}y^\varepsilon_{xx}(0,\cdot)- y_x^{1}(0,\cdot) \|_{L^2(0,T)} \leq C\varepsilon^{3/4},\\ 
& 
\| \varepsilon^{1/2}y^\varepsilon_{xx}(1,\cdot) +y_x^{1}(1,\cdot) \|_{L^2(0,T)} \leq C\varepsilon^{3/4}.
\end{aligned}
\end{equation}
\end{proposition}

\begin{remark}
The compatibility conditions in (\ref{eq_compo1}) come from the coupling between $y^0$ and $y^1$ in the boundary condition when solving (\ref{eq_sis11}). As $y^1\in C(0,T;H^4(0,1))$, we have that $y^1_{xxx}(0,t)\in C[0,T]$, and $y^0$ is even smoother. Then, the wave equation, satisfied by both $y^0$ and $y^1$ in the domain $Q_T$, has a trace at the boundary of $Q_T$ in this gives some constraints on the boundary conditions. Then,  
\begin{eqnarray*}
&&y^1(0,t)=-y^0_x(0,t), \quad t\in[0,T],\\
&&y^1_{xxx}(0,t)=y^1_{xtt}(0,t)=-y^0_{tt}(0,t)=-y^0_{xx}(0,t), \quad t\in[0,T],
\end{eqnarray*}
that must be true in particular at $t=0$. This provides the compatibility conditions in (\ref{eq_compo1}) for the initial data. An analogous argument at $x=1$ provides the rest of the compatibility conditions. 

Further compatibility conditions appear below for the second order term (see \ref{eq_compo1_3}) since $y^0$ and $y^1$ are assumed to be smoother and it appears a new function $y^2$ coupled with $y^1$ through the boundary condition in system (\ref{eq_sis11_2}).  

{\color{black} Note that, without these compatibility conditions on the initial data, solutions will have a new boundary layer along the characteristics starting at $x=0,1$. This requires an specific asymptotic analysis with a different ansatz (see \cite{AmiratMunch} for an example in this situation).}
\end{remark}

\begin{proposition} \label{pr_3}
Assume that $(y_0^{\varepsilon}, y_1^{\varepsilon}) \in X_1=H^2_0\times L^2$ and  consider, for $j=0,1,2$, $(y_0^{j}, y_1^{j}) \in Z_{6-j}\subset H^{6-j}\times H^{5-j}$, $ v^\varepsilon \in H^2_0(0,T)$ and $v^j \in H^{6-j}_0(0,T)$  that satisfies the following compatibility conditions
\begin{equation} \label{eq_compo1_3}
\begin{aligned}
& -y^0_{0,x} (0)=y_{0}^{1} (0), \quad y_{0,x}^{0} (1)=y_{0}^{1} (1), \\ 
& -y_{0,xxx}^{0} (0)=y_{0,xx}^{1} (0), \quad y_{0,xxx}^{0} (1)=y_{0,xx}^{1} (1), \\ 
& -y_{0,xxxxx}^{0} (0)=y_{0,xxxx}^{1} (0), \quad y_{0,xxxxx}^{0} (1)=y_{0,xxxx}^{1} (1),\\
& -y^1_{0,x} (0)=y_{0}^{2} (0), \quad y_{0,x}^{1} (1)=y_{0}^{2} (1), \\ 
& -y_{0,xxx}^{1} (0)=y_{0,xx}^{2} (0), \quad y_{0,xxx}^{1} (1)=y_{0,xx}^{2} (1).
\end{aligned}
\end{equation}

Let $y^\varepsilon$ be the solution of (\ref{eq_ep_gen}) and $y^0, \; y^1, \; y^2$ be the solutions of (\ref{eq_sis0}), (\ref{eq_sis11}) and (\ref{eq_sis11_2}) respectively. Define, 
\begin{eqnarray} \nonumber 
&&y^{\varepsilon,2}(x,t)= \sum_{j=0}^2 \varepsilon^{j/2} \left[ y^j(x,t)-
y^{j}(0,t) e^{-x/\varepsilon^{1/2}} -y^{j}(1,t) e^{-(1-x)/\varepsilon^{1/2}} \right] ,
\label{eq_app1_3}
\end{eqnarray}
and assume that,
\begin{equation} \label{eq_estdi2}
\begin{aligned}
& \| (y_0^{\varepsilon}, y_1^{\varepsilon})-(y_0^{\eps,2}, y_1^{\eps,2}) \|_{H^1\times L^2} +\varepsilon^{1/2} \| y_{0,xx}^{\varepsilon}-y_{0,xx}^{\eps,2}\|_{L^2(0,1)}\leq C\varepsilon^{3/2},\\ 
&\left\| \varepsilon^{1/2}v^\varepsilon - \sum_{j=0}^2 \varepsilon^{j/2} v^j\right\|_{H^2(0,T)} \leq C\varepsilon^{3/2}.
\end{aligned}
\end{equation}

Then, there exists a constant $C>0$, independent of $\varepsilon$, such that 
\begin{equation} \label{est_po1_3}
\begin{aligned}
& \| y^\varepsilon-y^{\varepsilon,2} \|_{L^\infty(0,T;H^1_0)} + \| y^\varepsilon_t-y^{\varepsilon,2}_t \|_{L^\infty(0,T;L^2)} + \varepsilon^{1/2} \| y^\varepsilon_{xx}- y^{\varepsilon,2}_{xx} \|_{L^\infty(0,T;L^2)}\leq C  \varepsilon^{5/4},\\
& \| \varepsilon^{1/2}y^\varepsilon_{xx}(0,\cdot)- y_x^{2}(0,\cdot) \|_{L^2(0,T)} \leq C\varepsilon^{5/4},\\ 
&\| \varepsilon^{1/2}y^\varepsilon_{xx}(1,\cdot) +y_x^{2}(1,\cdot) \|_{L^2(0,T)} \leq C\varepsilon^{5/4}.
\end{aligned}
\end{equation}
\end{proposition}

\section{Asymptotic expansion of the solution of the optimality system}\label{optimality_system}

We are now in position to determine the asymptotic expansion of the optimality system (\ref{eq_optsyseps1})-(\ref{eq_optsyseps2}) associated to the control of minimal $L^2(0,T)$-norm recalled in Section \ref{weightedcontrol}.
This system, through the equation $y^{\eps}(1,\cdot)=\eta \eps \ph^{\eps}_{xx}(1,\cdot)$ in $(0,T)$ makes the link between the direct and adjoint solution. Furthermore, as discussed in Section \ref{weightedcontrol}, the introduction of the weight time function $\eta$ makes the control $v^{\eps}$ more regular than $L^2(0,T)$ when associated to smooth initial $(y_0,y_1)$. This extra regularity allows to rigorously justify the asymptotic analysis. Along this section assume that both $(y^0,y^1)$ and $(\varphi^\varepsilon_0,\varphi^\varepsilon_1)$ are sufficiently smooth and satisfy the hypotheses of Proposition \ref{pr_3}. This will be justified in the next section. 

Note that in particular we assume that 
\begin{eqnarray} \nonumber 
&&\varepsilon^{1/2}\sum_{j=0}^2 \varepsilon^{j/2}v^j+\mathcal{O}(\varepsilon^{3/2})=\varepsilon^{1/2}\eta(t)\varphi^\varepsilon_{xx}(1,t)= -\eta(t)\sum_{j=0}^2 \varepsilon^{j/2}\varphi^j_x(1,t)+\mathcal{O}(\varepsilon^{3/2}), 
\end{eqnarray}
that we substitute in (\ref{eq_optsyseps1})-(\ref{eq_optsyseps2}). 
Then, we have 
\begin{equation} \label{eq_controles}
v^j(t)=-\eta (t) \varphi^j_x(1,t), \quad j=0,1,2.
\end{equation}

The aim of this section is to characterize $\varphi^j$ in terms of suitable optimality system associated to control problems for the wave equation. 

\subsection{Characterization of $v^0$}

According to Proposition \ref{pr_1} we can write the zero order asymptotics for the solution of (\ref{eq_optsyseps1})-(\ref{eq_optsyseps2})
\begin{eqnarray*}
&& \varphi^\varepsilon=\varphi^0+\mathcal{O}(\varepsilon^{1/2}), \\
&& \varphi_{xx}^\varepsilon(1,t)=-\varphi_x^0(1,t)+\mathcal{O}(\varepsilon^{1/2}), \\
&& y^\varepsilon(x,t)=y^0(x,t)-\eta (t)\varphi^0_x(1,t)e^{-(1-x)/\varepsilon^{1/2}}+\mathcal{O}(\varepsilon^{1/2}),
\end{eqnarray*} 
where $\varphi^0$ and $y^0$ are the solutions of the system:
\begin{equation}
 \label{eq_optsyswe1b}
 \left\{
 \begin{aligned}
& \varphi_{tt}^0  -  \varphi_{xx}^0 =0, & \text{ in } Q_T, \\
& \varphi^0 (0,\cdot)=\varphi_x^0 (0,\cdot)=0, & \text{ in } (0,T),\\
& \varphi^0(x,0)=\varphi_0^{0}, \quad \varphi_t^0(x,0)=\varphi_1^{0}, &\text{ in } \Omega,
\end{aligned}
\right. 
\end{equation}
\begin{equation}
  \label{eq_optsyswe2b}
\left\{
\begin{aligned}
& y_{tt}^0  -  y_{xx}^0 =0, & \text{ in } Q_T, \\
& y^0 (0,\cdot)=0, \quad  y^0 (1,\cdot)=-\eta\,\varphi_{x}^0(1,\cdot), & \text{ in } (0,T)\\
& y^0(\cdot,0)=y_0, \quad y_t^0(\cdot,0)=y_1, & \text{ in } \Omega, \\
& y(\cdot,T)= y_t(\cdot,T)=0, & \text{ in } \Omega.
\end{aligned}
\right.
\end{equation}

Therefore, $(\varphi_0^{0},\varphi_1^{0})$ satisfies the optimality system (\ref{eq_optsyswe1b})-(\ref{eq_optsyswe2b}). If we change the sign to $\varphi^0$ we obtain the optimality system (\ref{eq_optsyswe1}) associated to the wave equation. Therefore, $v^0=-\eta(t)\varphi_x^0(1,t)$ where $\eta(t)\varphi_x^0(1,t)$ is the control for the limit wave equation.

\subsection{Characterization of $v^1$}
 
{\color{black} According to Proposition \ref{pr_2} we can write the first order asymptotics for the solution of (\ref{eq_optsyseps1})-(\ref{eq_optsyseps2})
\begin{eqnarray*}
&& \varphi^\varepsilon=\varphi^0+\varepsilon^{1/2}\varphi^1 + \mathcal{O}(\varepsilon), \\
&& \varphi_{xx}^\varepsilon(1,t)=-\varphi_x^0(1,t)- \varepsilon^{1/2} \varphi_x^1(1,t)+\mathcal{O}(\varepsilon), \\
&& y^\varepsilon(x,t)=y^0(x,t)-\eta (t)\varphi^0_x(1,t)e^{-(1-x)/\varepsilon^{1/2}} \\
&& \quad
+ \varepsilon^{1/2}\left[y^1(x,t)-y_x^0(0,t)e^{-x/\varepsilon^{1/2}}+\left(y_x^0(1,t)-\eta (t)\varphi^1_x(1,t)\right)e^{-(1-x)/\varepsilon^{1/2}} \right]+\mathcal{O}(\varepsilon).
\end{eqnarray*} }

The new functions $\varphi^1$ and $y^1$ are the solutions of the system:
\begin{eqnarray} \label{eq_optsyswe1c}
&& \left\{
\begin{array}{ll}
\varphi_{tt}^1  -  \varphi_{xx}^1 =0, & \mbox{in $Q_T$} \\
\varphi^1 (0,t)=-\varphi_x^0 (0,t), \quad \varphi^1(1,t)=\varphi^0_x(1,t), & t\in (0,T)\\
\varphi^1(x,0)=\varphi_0^{1}, \quad \varphi_t(x,0)=\varphi_1^{1}, & x\in\Omega,
\end{array}
\right. \\  \label{eq_optsyswe2c}
&& \left\{
\begin{array}{ll}
y_{tt}^1  -  y_{xx}^1 =0, & \mbox{in $Q_T$} \\
y^1 (0,\cdot)=-y_x^0(0,\cdot), \quad  y^1 (1,\cdot)=y^0_x(1,\cdot)-\eta(t)\varphi^1_x(1,\cdot),& t\in (0,T)\\
y^1(\cdot,0)= y_t^1(\cdot,0)=0, & x\in\Omega, \\
y^1(\cdot,T)= y_t^1(\cdot,T)=0, & x\in\Omega,
\end{array}
\right.
\end{eqnarray}
Here we decompose  $\ph^1=\ph+\ph^a$ with 
\begin{equation}
\label{eq:phia}
\left\{
\begin{aligned}
& \ph_{tt}- \ph_{xx}=0, &\text{ in }  \, Q_T, \\
& \ph(0,\cdot)= \ph(1,\cdot)=0 &\text{ in }  \, (0,T),   \\
&  (\ph(\cdot,0),\ph_t(\cdot,0))=(\ph^1_0,\ph^1_1), &\text{ in }  \, \Omega.
 \end{aligned}
\right.\qquad 
\left\{
\begin{aligned}
& \ph^{a,1}_{tt}- \ph_{xx}^{a,1}=0, &\text{ in }  \, Q_T, \\
& \ph^{a,1}(0,\cdot)=-\ph_x^0(0,\cdot), &\text{ in }  \, (0,T), \\
& \ph^{a,1}(1,\cdot)=\ph_x^0(1,\cdot), &\text{ in }  \, (0,T),   \\
&  (\ph^{a,1}(\cdot,0),\ph^{a,1}_t(\cdot,0))=(0,0), &\text{ in }  \, \Omega.
 \end{aligned}
\right.
\end{equation}
so that $\ph_x^1(1,\cdot)=\ph_x(1,\cdot)+\ph_x^{a,1}(1,\cdot)$. The system (\ref{eq_optsyswe1c}) then becomes 
\begin{eqnarray} \label{eq_optsyswe1cc}
&& \left\{
\begin{array}{ll}
\varphi_{tt}  -  \varphi_{xx} =0, & \mbox{in $Q_T$} \\
\varphi (0,\cdot)=\varphi(1,\cdot)=0, & t\in (0,T)\\
\varphi(\cdot,0)=\varphi_0^{1}, \quad \varphi_t(\cdot,0)=\varphi_1^{1}, & \in\Omega,
\end{array}
\right. \\  \label{eq_optsyswe2cc}
&& \left\{
\begin{array}{ll}
y_{tt}^1  -  y_{xx}^1 =0, & \mbox{in $Q_T$} \\
y^1 (0,\cdot)=-y_x^0(0,\cdot), \quad  y^1 (1,\cdot)=y^0_x(1,\cdot)-\eta(t)\varphi_x(1,\cdot)-\eta \ph^{a,1}_x(1,t),& t\in (0,T)\\
y^1(\cdot,0)= y_t^1(\cdot,0)=0, & x\in\Omega, \\
y^1(\cdot,T)= y_t^1(\cdot,T)=0, & x\in\Omega,
\end{array}
\right.
\end{eqnarray}

Therefore, $(\varphi_0^{1},\varphi_1^{1})$ satisfies the optimality system (\ref{eq_optsyswe1cc})-(\ref{eq_optsyswe2cc}). If we change the sign to $\varphi$ we obtain the optimality system (\ref{eq_optsyswe1}) associated to the wave equation (\ref{eq_optsyswe2cc}) with nonhomogenous boundary conditions. Therefore, $v^1=-\eta(t)\varphi_x(1,t)$ where $\eta(t)\varphi_x(1,t)$ is the control for this wave equation.

{\color{black} Note that this control $v^1=-\eta(t)\varphi_x(1,t)$ can be also characterized in terms of a boundary control for the wave equation with homogeneous boundary conditions and a particular initial data. In fact, if we define $(g^1_0,g^1_1)=(-\hat y^1(x,0), -\hat y^1_t(x,0))$ where $\hat y^1$ is the solution of the backwards system
$$
\left\{
\begin{array}{ll}
\hat y_{tt}^1  -  \hat y_{xx}^1 =0, & \mbox{in $Q_T$} \\
\hat y^1 (0,\cdot)=- y_x^0(0,\cdot), \quad  \hat y^1 (1,\cdot)= y^0_x(1,\cdot)-\eta \ph^{a,1}_x(1,t),& t\in (0,T)\\
\hat y^1(\cdot,T)= \hat y_t^1(\cdot,T)=0, & x\in\Omega,
\end{array}
\right.
$$
then, by linearity, we easily see that $v^1$ is the boundary control of the wave equation (\ref{eq:sys_y}) associated to the initial data $(g_0^1,g_1^1)$. 

\begin{remark} \label{re:regg} Note that the regularity of $v^1$ is one degree less than the one of $v^0$. In fact, $y^1$ has as boundary conditions the normal derivative of $y^0$ and therefore, one degree less of regularity than $y^0$. This is translated into the regularity of the controls for $y^1$. 
\end{remark} 
}

\subsection{Characterization of $v^2$}

The second order asymptotics for $y^\varepsilon$ and $\varphi^\varepsilon$ makes to appear the new functions $\varphi^2$ and $y^2$ which are the solutions of the system:
\begin{eqnarray} \label{eq_optsyswe1d}
&& \left\{
\begin{array}{ll}
\varphi_{tt}^2  -  \varphi_{xx}^2 =-\ph_{xxxx}^0, & \mbox{in $Q_T$} \\
\varphi^2 (0,t)=-\varphi_x^1 (0,t), \quad \varphi^2(1,t)=\varphi^1_x(1,t), & t\in (0,T)\\
\varphi^2(x,0)=\varphi_0^{2}, \quad \varphi_t^2(x,0)=\varphi_1^{2}, & x\in\Omega,
\end{array}
\right. \\  \label{eq_optsyswe2d}
&& \left\{
\begin{array}{ll}
y_{tt}^2  -  y_{xx}^2 =-y^0_{xxxx}, & \mbox{in $Q_T$} \\
y^2 (0,\cdot)=-y_x^1(0,\cdot), \quad  y^2 (1,\cdot)=y^1_x(1,\cdot)+\frac12 y^0_{tt}(1,\cdot)-\eta(t)\varphi^2_x(1,\cdot),& t\in (0,T)\\
y^2(\cdot,0)= y_t^2(\cdot,0)=0, & x\in\Omega, \\
y^2(\cdot,T)= y_t^2(\cdot,T)=0, & x\in\Omega,
\end{array}
\right.
\end{eqnarray}
Here we decompose  $\ph^2=\ph+\ph^{a,2}$ with 
\begin{equation}
\label{eq:phia2}
\left\{
\begin{aligned}
& \ph_{tt}- \ph_{xx}=0, &\text{ in }  \, Q_T, \\
& \ph(0,\cdot)=0, \ph(1,\cdot)=0 &\text{ in }  \, (0,T),   \\
&  (\ph(\cdot,0),\ph_t(\cdot,0))=(\ph^2_0,\ph^2_1), &\text{ in }  \, \Omega.
 \end{aligned}
\right.\qquad 
\left\{
\begin{aligned}
& \ph^{a,2}_{tt}- \ph_{xx}^{a,2}=-\ph_{xxxx}^0, &\text{ in }  \, Q_T, \\
& \ph^{a,2}(0,\cdot)=-\ph_x^1(0,\cdot), \ph^{a,2}(1,\cdot)=\ph_x^1(1,\cdot) &\text{ in }  \, (0,T),   \\
&  (\ph^{a,2}(\cdot,0),\ph^{a,2}_t(\cdot,0))=(0,0), &\text{ in }  \, \Omega.
 \end{aligned}
\right.
\end{equation}
so that $\ph_x^2(1,\cdot)=\ph_x(1,\cdot)+\ph_x^{a,2}(1,\cdot)$. The system (\ref{eq_optsyswe1d})-(\ref{eq_optsyswe2d}) then becomes 
\begin{eqnarray} \label{eq_optsyswe1dd}
&& \left\{
\begin{array}{ll}
\varphi_{tt}  -  \varphi_{xx} =0, & \mbox{in $Q_T$} \\
\varphi (0,\cdot)=\varphi(1,\cdot)=0, & t\in (0,T)\\
\varphi(\cdot,0)=\varphi_0^{1}, \quad \varphi_t(\cdot,0)=\varphi_1^{1}, & \in\Omega,
\end{array}
\right. \\  \label{eq_optsyswe2dd}
&& \left\{
\begin{array}{ll}
y_{tt}^2  -  y_{xx}^2 =-y^0_{xxxx}, & \mbox{in $Q_T$} \\
y^2 (0,\cdot)=-y_x^1(0,\cdot), \\  
y^2 (1,\cdot)=y^1_x(1,\cdot)+\frac12 y^0_{tt}(1,\cdot)-\eta(t)\varphi_x(1,\cdot)-\eta \ph^{a,2}_x(1,t),& t\in (0,T)\\
y^1(\cdot,0)= y_t^1(\cdot,0)=0, & x\in\Omega, \\
y^1(\cdot,T)= y_t^1(\cdot,T)=0, & x\in\Omega.
\end{array}
\right.
\end{eqnarray}

Therefore, $(\varphi_0^{2},\varphi_1^{2})$ satisfies the optimality system (\ref{eq_optsyswe1cc})-(\ref{eq_optsyswe2cc}). If we change the sign to $\varphi$ we obtain the optimality system (\ref{eq_optsyswe1}) associated to the wave equation (\ref{eq_optsyswe2cc}) with nonhomogenous boundary conditions. Therefore, $v^2=-\eta(t)\varphi_x(1,t)$ where $\eta(t)\varphi_x(1,t)$ is the control for this wave equation.

{\color{black} As in the previous case, we can also write  $v^2=-\eta(t)\varphi_x(1,t)$ as the boundary control for a wave equation with homogeneous boundary condition and suitable initial data. On the other hand, as pointed out in Remark \ref{re:regg} for $v^1$, the control $v^2$ will have one less degree of regularity than $v^1$. }

\section{Convergence of the controls}\label{convergence_control}

In this section we prove the following result
\begin{theorem}\label{estimate_errorcontrol}
Let $0\leq n\leq 2$. Assume that $(y_0, y_1)\in Z_{4+n}$ and satisfies the compatibility condition in (\ref{eq_compo1}) (for $n=1$) or (\ref{eq_compo1_3}) (for $n=2$).
Consider $v^j$, $0\leq j\leq n$, the controls obtained in the previous section. Let $\eps>0$ and $v^{\eps}$ be the control of minimal $L^2$-weighted norm for (\ref{eq:sys_yeps}) associated to the data $(y_0,y_1)$. Then, there exists a constant $C>0$ such that
\begin{equation}\nonumber
\biggl\| \varepsilon^{1/2}v^\varepsilon - \sum_{j=0}^n \varepsilon^{j/2}v^j \biggr\|_{L^2(0,T)} \leq C\varepsilon^{n/2+1/4}.
\end{equation} 
\end{theorem}

\par\noindent
{\sc  Proof-} We follow the notation of the previous section. The strategy is to show that $\varepsilon^{1/2}v^\varepsilon - \sum_{j=0}^2 \varepsilon^{j/2}v^j$ is a control of minimal $L^2$-morm which drives to rest a solution of (\ref{eq:sys_yeps})  associated to a vanishing initial condition as $\eps$ goes to zero. The result then follows from the continuous dependance of control of minimal $L^2$-norm with respect to the initial condition to be controlled. We divide the proof in two steps. 

\bigskip
\par\noindent
STEP 1. Let $(\varphi^{j}_0,\varphi^{j}_1)$, $0\leq j\leq 2$, be the initial data that provides the controls $v^j$ in the previous section. Note that the initial data $(y_0, y_1)\in Z_{4+n}$ and then we have 
\begin{equation} \label{eq_regfi}
(\varphi^{j}_0,\varphi^{j}_1)\in Z_{5+n-j}, \quad j\leq n. 
\end{equation}
In fact, $(\varphi^{0}_0,\varphi^{0}_1)$ is characterized by the optimality system (\ref{eq_optsyswe1b})-(\ref{eq_optsyswe2b}) which corresponds to a control for the wave equation with initial data $(y_0, y_1)\in Z_{4+n}$. Therefore, as stated in Theorem \ref{th_EZwe}, $(\varphi^{0}_0,\varphi^{0}_1)\in Z_{5+n}$. On the other hand,  $(\varphi^{1}_0,\varphi^{1}_1)$ is characterized by the optimality system (\ref{eq_optsyswe1c})-(\ref{eq_optsyswe2c}) which corresponds to a control for a solution of the wave equation with one degree less of regularity {\color{black} (see Remark \ref{re:regg} above)}. Therefore, $(\varphi^{1}_0,\varphi^{1}_1) \in Z_{4+n}$, and so on.  

 Consider the following initial data for the adjoint system (\ref{eq:sys_pheps1d}), 
\begin{eqnarray} \label{eq_psi_01} 
\psi^{\varepsilon}_0(x)&=& \sum_{j=0}^n \varepsilon^{j/2} \left( 
\varphi^{j}_0(x) -\varphi^{j}_0(0) e^{-x/\varepsilon^{1/2}} 
- \varphi^{j}_0(1) e^{-(1-x)/\varepsilon^{1/2}} \right), \\ \label{eq_psi_02} 
\psi^{\varepsilon}_1(x)&=& \sum_{j=0}^n \varepsilon^{j/2} \left( 
\varphi^{j}_1(x) -\varphi^{j}_1(0) e^{-x/\varepsilon^{1/2}} 
- \varphi^{j}_1(1) e^{-(1-x)/\varepsilon^{1/2}} \right).
\end{eqnarray}

In this step we prove that this final data is associated to a suitable control that we characterize below.

The first difficulty is that  $(\psi^{\varepsilon}_0,\psi^{\varepsilon}_1)$, as defined in (\ref{eq_psi_01})-(\ref{eq_psi_02}) is not in $ H^2_0\times L^2$, in general, due to the fact that $\psi^{\varepsilon}_0$ does not satisfy the homogeneous boundary condition. This is related to the composite method we employed to define the approximations of $y^{\eps}$ and $\ph^\eps$. In fact, 
\begin{equation}
\begin{aligned}
\psi^{\varepsilon}_0(0)&=-\sum_{j=0}^n \varepsilon^{j/2}  
\varphi^{j}_0(1) e^{-1/\varepsilon^{1/2}} =O(e^{-1/\varepsilon^{1/2}}) ,\\
\psi^{\varepsilon}_{0,x}(0)&=\sum_{j=0}^n \varepsilon^{j/2}\left( 
\varphi^{j}_{0,x}(0) +\varepsilon^{-1/2}\varphi^{j}_0(0) - \varepsilon^{-1/2} \varphi^{j}_0(1) e^{-1/\varepsilon^{1/2}} \right)\\
&=\varepsilon^{-1/2}\ph^0_0(0)+ \sum_{j=0}^{n-1} \varepsilon^{j/2} \left(\varphi^{j}_{0,x}(0) + \varphi^{j+1}_0(0) \right) + \varepsilon^{n/2}\varphi^{n}_{0,x}(0)+\mathcal{O}(e^{-1/\varepsilon^{1/2}}) \\
&= \varepsilon^{n/2}\varphi^{n}_{0,x}(0)+\mathcal{O}(e^{-1/\varepsilon^{1/2}}).
\end{aligned}
\end{equation}
An analogous situation appears at the boundary at $x=1$. 

To overcome this difficulty we correct the right hand side in (\ref{eq_psi_01}) by a function $R^{T,\varepsilon}$ in such a way that 
$$
\psi^{T,\varepsilon}_0(x)= \sum_{j=0}^n \varepsilon^{j/2} \left( 
\varphi^{T,j}_0(x) -\varphi^{T,j}_0(0) e^{-x/\varepsilon^{1/2}} 
- \varphi^{T,j}_0(1) e^{-(1-x)/\varepsilon^{1/2}} \right) +R^{T,\varepsilon}(x),
$$
satisfies the boundary conditions and 
$$
\| R^{T,\varepsilon}\|_{H^1(0,1)} + \varepsilon^{1/2} \| R^{T,\varepsilon} \|_{H^2(0,1)}\leq C \varepsilon^{n/2+1/4}.
$$
This can be achieved with the function 
$$
R^{T,\varepsilon}(x)=\varepsilon^{n/2}\varphi^{n}_{0,x}(0)e^{-x/\varepsilon^{1/2}}
-\varepsilon^{n/2}\varphi^{n}_{0,x}(1)e^{-(1-x)/\varepsilon^{1/2}}+P(x)e^{-1/\varepsilon^{1/2}},
$$
where $P$ is a suitable polynomial that takes into account all the terms with the factor $e^{-1/\varepsilon^{1/2}}$.

Note that our choice of $(\psi^{T,\varepsilon}_0, \psi^{T,\varepsilon}_1)$ satisfies the hypothesis of Propositions \ref{pr_2} or \ref{pr_3} (depending on $n$). For example, if $n=2$ the regularity asumptions  and compatibility conditions hold in view of (\ref{eq_regfi}). Then,  the solution $\psi$ of the adjoint system satisfies, 
\begin{equation} \label{eq_psid}
\psi^\varepsilon(x,t) =\sum_{j=0}^n \varepsilon^{j/2} \left(\varphi^j(x,t)- \varphi^j(0,t)e^{-x/\varepsilon^{1/2}}-\varphi^j(1,t) e^{-(1-x)/\varepsilon^{1/2}}  \right) + \mathcal{O}(\varepsilon^{n/2+1/4})  ,
\end{equation}
and 
\begin{eqnarray} \label{eq_as_whumpsi}
&&\varepsilon^{1/2}\psi^\varepsilon_{xx}(1,t)= -\sum_{j=0}^n \varepsilon^{j/2}\varphi^j_x(1,t)+\mathcal{O}(\varepsilon^{n/2+1/4})=\sum_{j=0}^n \varepsilon^{j/2}v^j+\mathcal{O}(\varepsilon^{n/2+1/4}). 
\end{eqnarray}
Let us define the system 
\begin{equation} \label{eq_sis_1_pf}
\left\{ \begin{array}{ll}
z^\varepsilon_{tt} +\varepsilon z^\varepsilon_{xxxx}- z^\varepsilon_{xx} =0, & \mbox{in}\quad Q_T, \\
z^\varepsilon (0,t)=z^\varepsilon(1,t)=0,  & t\in(0,T) \\
z^\varepsilon_x(0,t)=0, \quad z^\varepsilon_x(1,t)=\eta(t)\psi_{xx}^\varepsilon(1,t),  & t\in(0,T) \\
z^\varepsilon(x,0)= y_0(x), \quad  z^\varepsilon_{t}(x,0)=y_1(x), & x\in(0,1),
\end{array} \right.
\end{equation}
and 
$(g^\varepsilon_0,g^\varepsilon_1)=(z^\varepsilon(\cdot,T),z^\varepsilon_t(\cdot,T))$. By Proposition \ref{pr_2} and our choice of $(\psi_0^{\varepsilon},\psi_1^{\varepsilon})$ we have 
$$
\varepsilon^{1/2} \| g_{0,xx}^0 \|_{L^2} +\|(g^\varepsilon_0,g^\varepsilon_1)\|_{H^1\times L^2} = \mathcal{O}(\varepsilon^{n/2+1/4}).
$$
Observe that $\eta(t)\psi_{xx}^\varepsilon(1,t)$ is a control for (\ref{eq_sis_1_pf}) that drives the initial data $(y_0,y_1)$ to $(g^\varepsilon_0,g^\varepsilon_1)$.

\bigskip
\par\noindent
STEP 2. Consider now the function $w^\varepsilon=y^\varepsilon - z^\varepsilon$ and $\zeta^\varepsilon = \varphi^\varepsilon-\psi^\varepsilon$. They satisfy the coupled system 
\begin{eqnarray} \label{eq_sis_1_pf2}
&& \left\{ \begin{array}{ll}
\zeta^\varepsilon_{tt} +\varepsilon \zeta^\varepsilon_{xxxx}- \zeta^\varepsilon_{xx} =0, & \mbox{in}\quad Q_T, \\
\zeta^\varepsilon (0,t)=\zeta^\varepsilon(1,t)=0,  & t\in(0,T), \\
\zeta^\varepsilon_x(0,t)=\zeta^\varepsilon_x(1,t)=0,  & t\in(0,T), \\
\zeta^\varepsilon(x,0)=\varphi^\varepsilon_0-\psi^\varepsilon_0, \quad  \zeta^\varepsilon_{t}(x,0)=\varphi^\varepsilon_1-\psi^\varepsilon_1, & x\in(0,1),
\end{array} \right. \\
&& \left\{ \begin{array}{ll}
w^\varepsilon_{tt} +\varepsilon w^\varepsilon_{xxxx}- w^\varepsilon_{xx} =0, & \mbox{in}\quad Q_T, \\
w^\varepsilon (0,t)=w^\varepsilon(1,t)=0,  & t\in(0,T), \\
w^\varepsilon_x(0,t)=0, \quad w^\varepsilon_x(1,t)=\eta(t)\zeta_{xx}^\varepsilon(1,t),  & t\in(0,T), \\
w^\varepsilon(x,0)= 0, \quad  w^\varepsilon_{t}(x,0)=0, & x\in(0,1), \\
w^\varepsilon(x,T)=-g_0^\varepsilon, \quad w^\varepsilon_t(x,T)=-g_1^\varepsilon .
\end{array} \right. 
\end{eqnarray} 
Note that this is the optimality system for the unique minimal weighted $L^2$-norm that drives the initial state $(0,0)$ to the final state $(-g_0^\varepsilon,-g_1^\varepsilon)$
Therefore, by estimate (\ref{eq_est_cc}) we obtain
\begin{eqnarray*}
&& \| \eta(t) \zeta^\varepsilon \|_{L^2(0,T)} =
\| v^\varepsilon-\eta(t)\psi_{xx}^\varepsilon(1,t) \|_{L^2(0,T)} \leq C \|(g^\varepsilon_0,g^\varepsilon_1)\|_{X_1} \\
&& \quad =\varepsilon^{1/2} \| g_{0,xx}^0 \|_{L^2} +\|(g^\varepsilon_0,g^\varepsilon_1)\|_{H^1\times L^2} = \mathcal{O}(\varepsilon^{n/2+1/4})
\end{eqnarray*}
which allows to conclude. $\hfill\Box$

\begin{remark}
Once proved the convergence of the controls stated in Theorem \ref{estimate_errorcontrol} one can easily state a convergence result for the controlled solutions in the energy space, thanks to Propositions \ref{pr_1}, \ref{pr_2} or \ref{pr_3} (depending on $n$). 
\end{remark}

\section{Numerical experiments}\label{experiments}

We illustrate our theoretical results with one numerical experiment. Precisely, we take $T=2.5$ and the initial condition $(y_0,y_1)=(\sin(2\pi x)^4,0)\in Z^6$. We consider a weight function $\eta$ defined as follows 
\begin{equation}\nonumber
\eta(t)=\big((1-e^{-40 t})(1-e^{-40(T-t)})\big)^3, \quad t\in [0,T].
\end{equation}
For any $\eps$ fixed, the control $v^\eps$ of minimal $L^2(0,T; \eta)$ weight norm for the system (\ref{eq:sys_yeps}) is computed by minimizing the conjugate functional $J_\eps^{\star}$ with respect to the initial condition of the adjoint state. This is performed using the Polak-Ribiere version of  conjugate gradient method. The iterative process is stopped when the sequence $\{y_k^\eps\}_{k>0}$ related to the minimizing sequence of $J_\eps^\star$ satisfies
\begin{equation}
\Vert y^{\eps}_{k_0}(\cdot,T),y_{k_0,t}^{\eps}(\cdot,T)\Vert_{H^2(0,1)\times L^2(0,1)}\leq 10^{-6}, \label{stop}
\end{equation}
for some $k_0=k_0(\eps)\in \mathbb{N}$.

A $C^1$-finite element approximation is used for the space variable and a centered finite difference scheme is used for the time variable. We refer to \cite{munch_arch2010} for the details in the similar context of the linear system (\ref{arch_essentialspectrum}) employing the general approach discussed in \cite{Glowinski_He_Lions,Munch_M2AN05}. A similar method approach is used to approximate the controls of Dirichlet type $v^0$,$v^1$ and $v^2$ based on \cite{cindea_munch_calcolo2015}.

\begin{figure}[http!]
\begin{center}
\includegraphics[width=9cm,height=7cm]{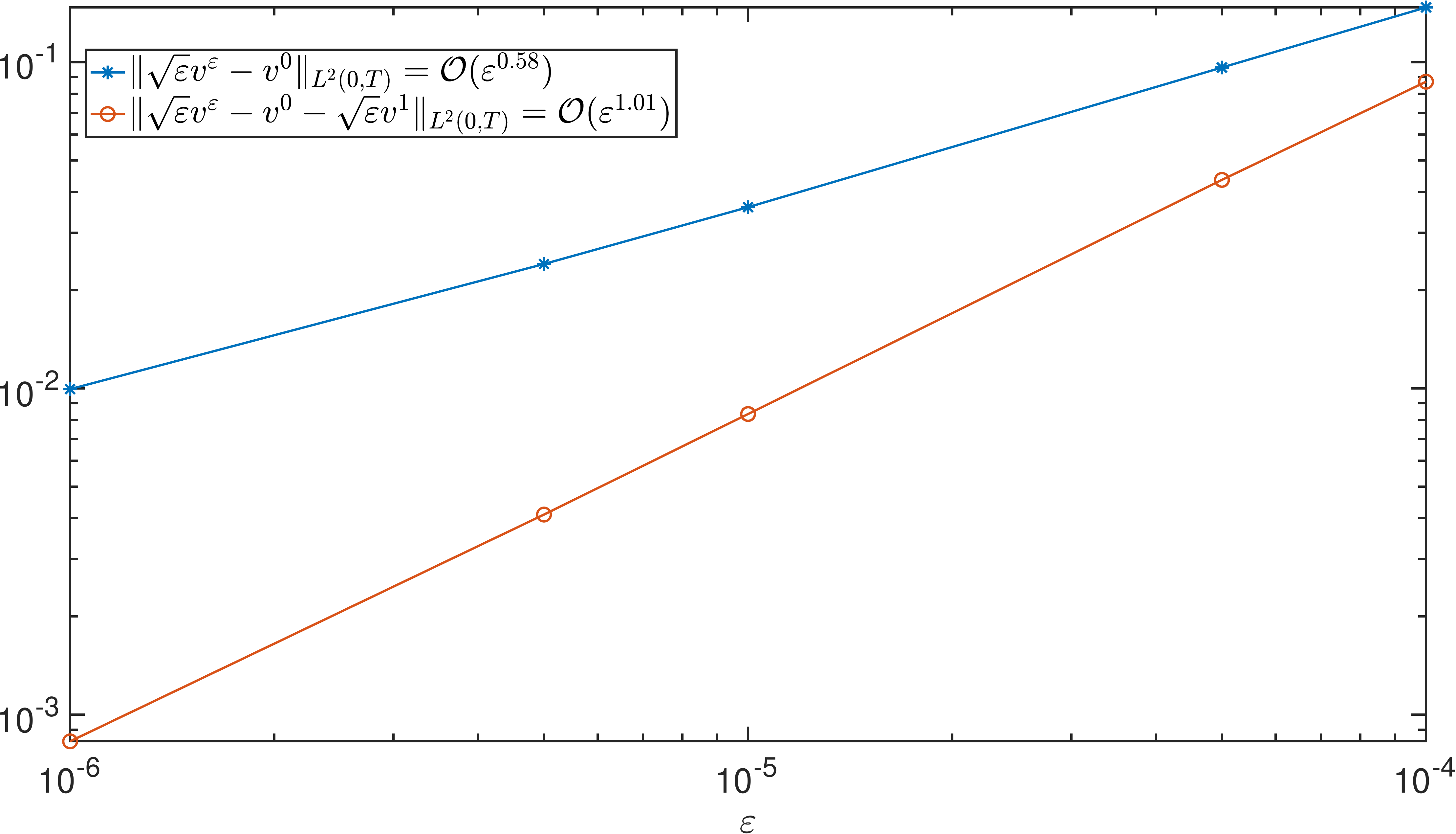}
\caption{Evolution of $\Vert \sqrt{\eps}v^{\eps}-v^0\Vert_{L^2(0,T)}$ and $\Vert \sqrt{\eps}v^{\eps}-v^0-\sqrt{\eps}v^1\Vert_{L^2(0,T)}$ with respect to $\eps$.}\label{fig:rate}
\end{center}
\end{figure}

Table \ref{tab:norms} collects the norm of $\sqrt{\eps}v^{\eps}$ for values of $\eps\in (10^{-6},10^{-1})$ and highlights the uniform bound property of on $\{\sqrt{\eps}\Vert v^{\eps}\Vert_{L^2(0,T)}\}_{\eps>0}$ according to the earlier results du to J.-L.Lions (Theorem \ref{theoremLions}). The table also emphasizes that the number of iterates $k_0(\eps)$ to achieve (\ref{stop}) increases as $\eps\to 0$ and traduces the loss of the uniform coercivity of the functional $J_\eps^\star$ with respect to the norm $(H^2\cap H_0^1)(0,1)\times L^2(0,1)$. The remaining part of the table gives the error $E^{\eps}_k=\Vert \sqrt{\eps}v^{\eps}-\sum_{n=0}^k \eps^{n/2}v^n\Vert_{L^2(0,T)}$: for $\eps$ small enough, we compute 
$$
\begin{aligned}
& E^\eps_0=\Vert \sqrt{\eps}v^{\eps}-v^0\Vert_{L^2(0,T)}=\mathcal{O}(\eps^{0.58}),\\
& E^\eps_1=\Vert \sqrt{\eps}v^{\eps}-v^0-\sqrt{\eps}v^1\Vert_{L^2(0,T)}=\mathcal{O}(\eps^{1.01}),\\
& E^\eps_2=\Vert \sqrt{\eps}v^{\eps}-v^0-\sqrt{\eps}v^1-\eps v^2\Vert_{L^2(0,T)}=\mathcal{O}(\eps^{1.36}),
\end{aligned}
$$
so that we observe slightly better rates than those given in Theorem \ref{estimate_errorcontrol}. We also refer to Figure \ref{fig:rate}. 

\begin{table}[http]
	\centering
		\begin{tabular}{|c|c||c|c|c|c|}
			\hline
			$\eps$  			 & $\sharp$ iterates & $\Vert \sqrt{\eps}v_{\eps}\Vert_{L^2(0,T)}$ &   $E^\eps_0$ & $E^\eps_1$ & $E^\eps_2$\tabularnewline
			\hline
			$10^{-1}$ 	        & 5	& $0.2625$	& $4.68\times 10^{-1}$ & $4.12\times 10^{-1}$ &  $3.1\times 10^{-1}$ \tabularnewline
			$10^{-2}$            & 11 & $0.2965$ &    $4.28\times 10^{-1}$ & $3.32\times 10^{-1}$ &  $2.1\times 10^{-1}$\tabularnewline
			$10^{-3}$ 	       & 24 & $0.3542$ &		 $3.61\times 10^{-1}$ &  $2.82\times 10^{-1}$ & $1.79\times 10^{-1}$	\tabularnewline
			$10^{-4}$ 	        & 51 & $0.3510$ &		 $1.47\times 10^{-1}$ & $8.71\times 10^{-2}$ & 	$6.21\times 10^{-2}$\tabularnewline
			$5\times 10^{-5}$ &	 90       & $0.3508$ &		 $9.29\times 10^{-2}$ & $4.35\times 10^{-2}$ & $2.01\times 10^{-2}$ \tabularnewline
			$10^{-5}$ 	        & 101 & $0.3499$ &		 $3.59\times 10^{-2}$& $8.34\times 10^{-3}$  & $2.37\times 10^{-3}$	\tabularnewline
			$5\times 10^{-6}$ &	171        & $0.3498$ &		 $2.40\times 10^{-2}$& $4.30\times 10^{-3}$  & $9.31\times 10^{-4}$	\tabularnewline
			$10^{-6}$ 	        & 203 & $0.3498$ &		 $9.95\times 10^{-3}$  &  $8.34\times 10^{-4}$ & $1.13\times 10^{-4}$\tabularnewline
                       			\hline
		\end{tabular}
	\caption{$L^2(0,T)$ norms of the $\sqrt{\eps}v^\eps$ and of the error $E^\eps_k=\Vert \sqrt{\eps}v^{\eps}-\sum_{n=0}^k \eps^{n/2}v^n\Vert_{L^2(0,T)}$, $k=0,1,2$,  with respect to $\eps\in (10^{-6},10^{-1})$. $C^0$ case - $\Vert v^0\Vert_{L^2(0,T)} \approx 0.349834$.}
	\label{tab:norms}
	\end{table}
	
\begin{figure}[http]
\begin{center}
\includegraphics[scale=0.4]{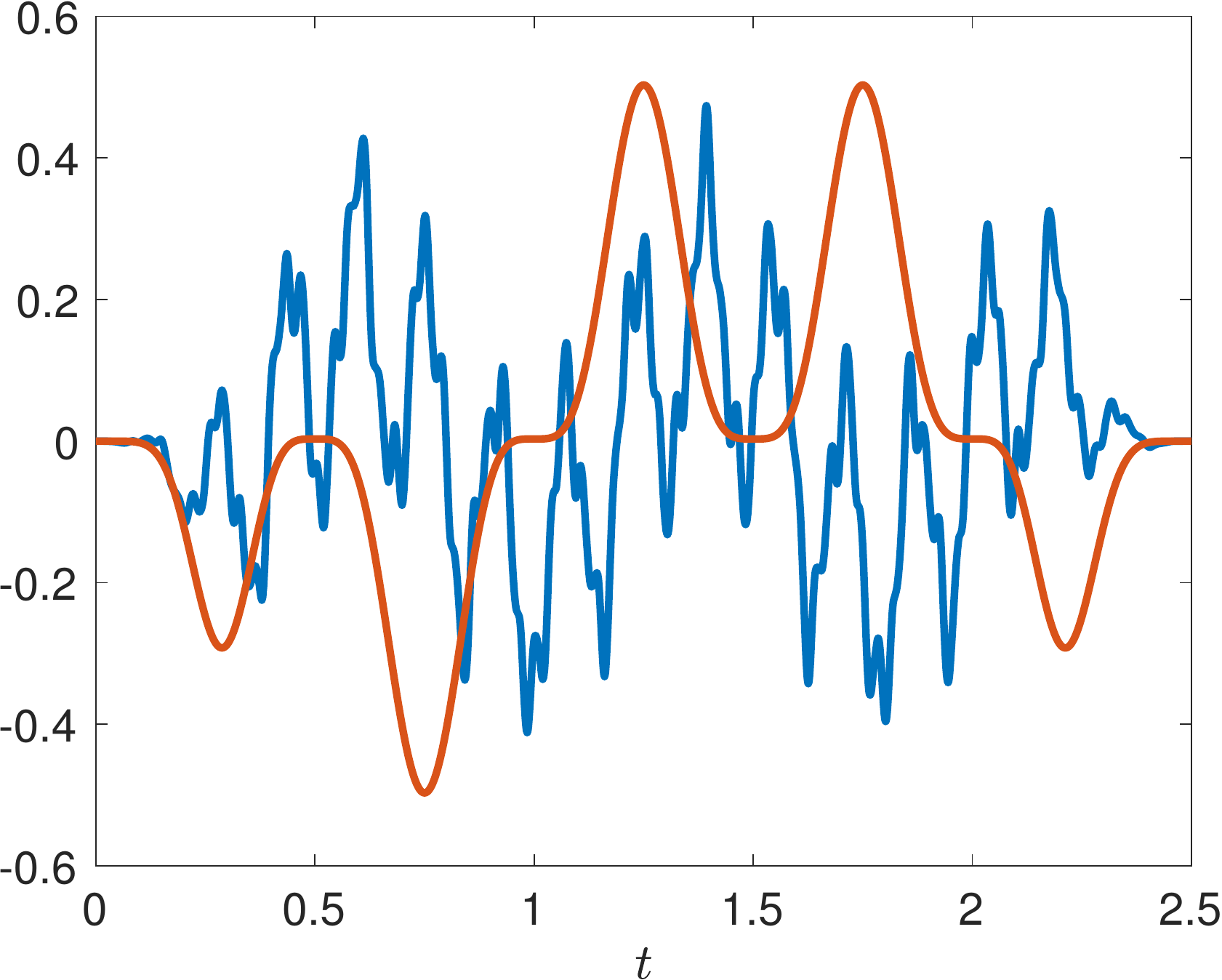}
\hspace*{0.3cm}
\includegraphics[scale=0.4]{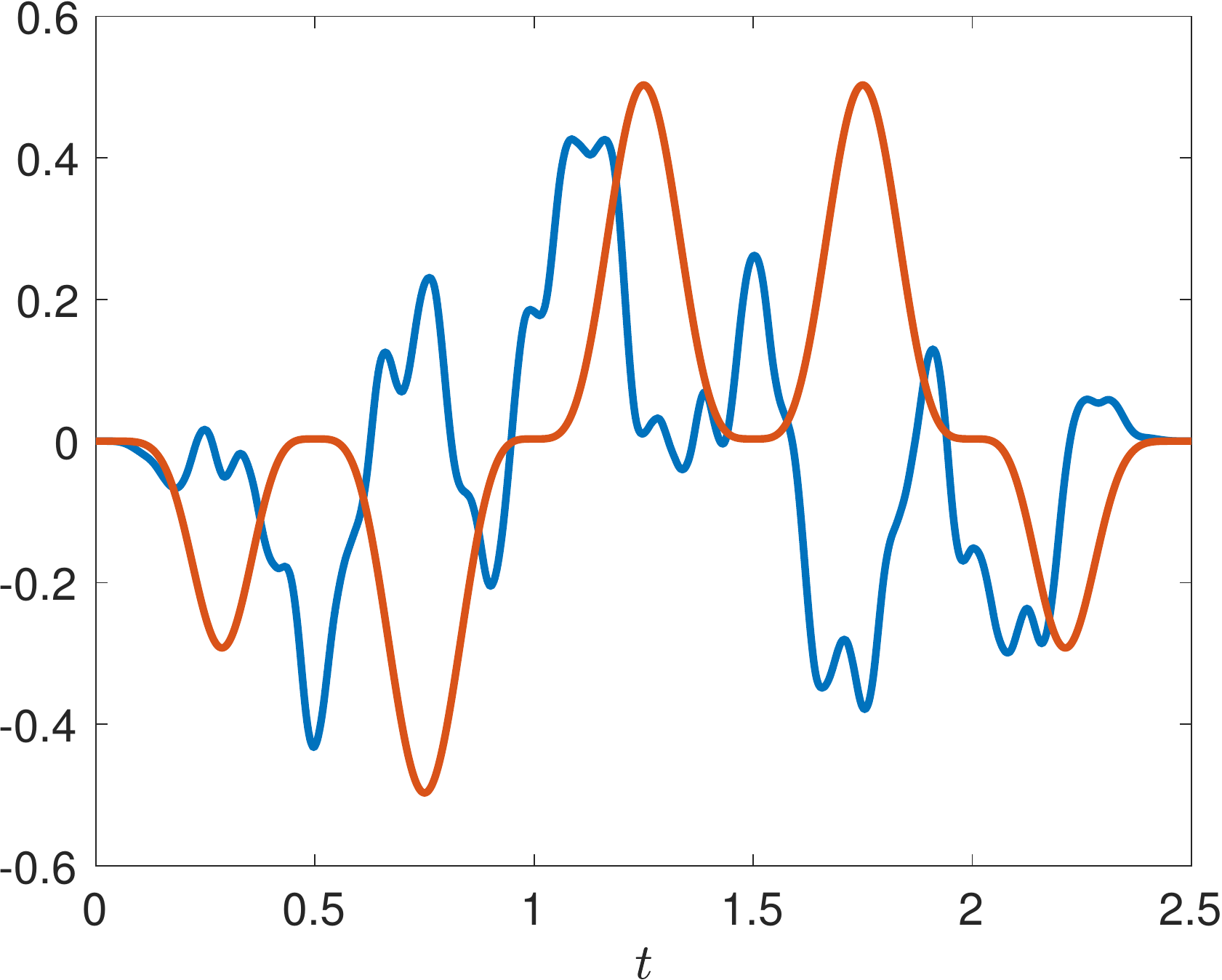}\\
\includegraphics[scale=0.4]{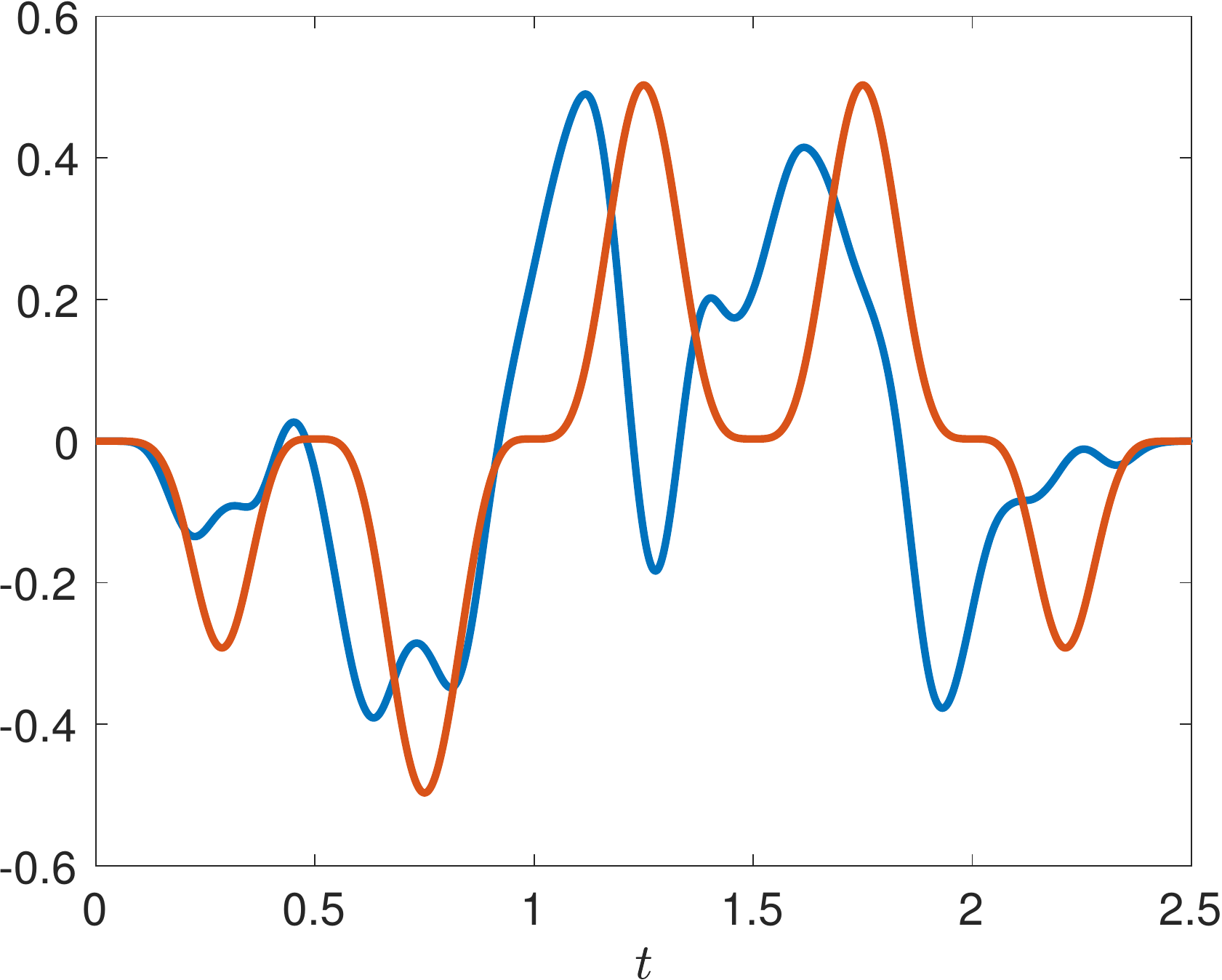}
\hspace*{0.3cm}
\includegraphics[scale=0.4]{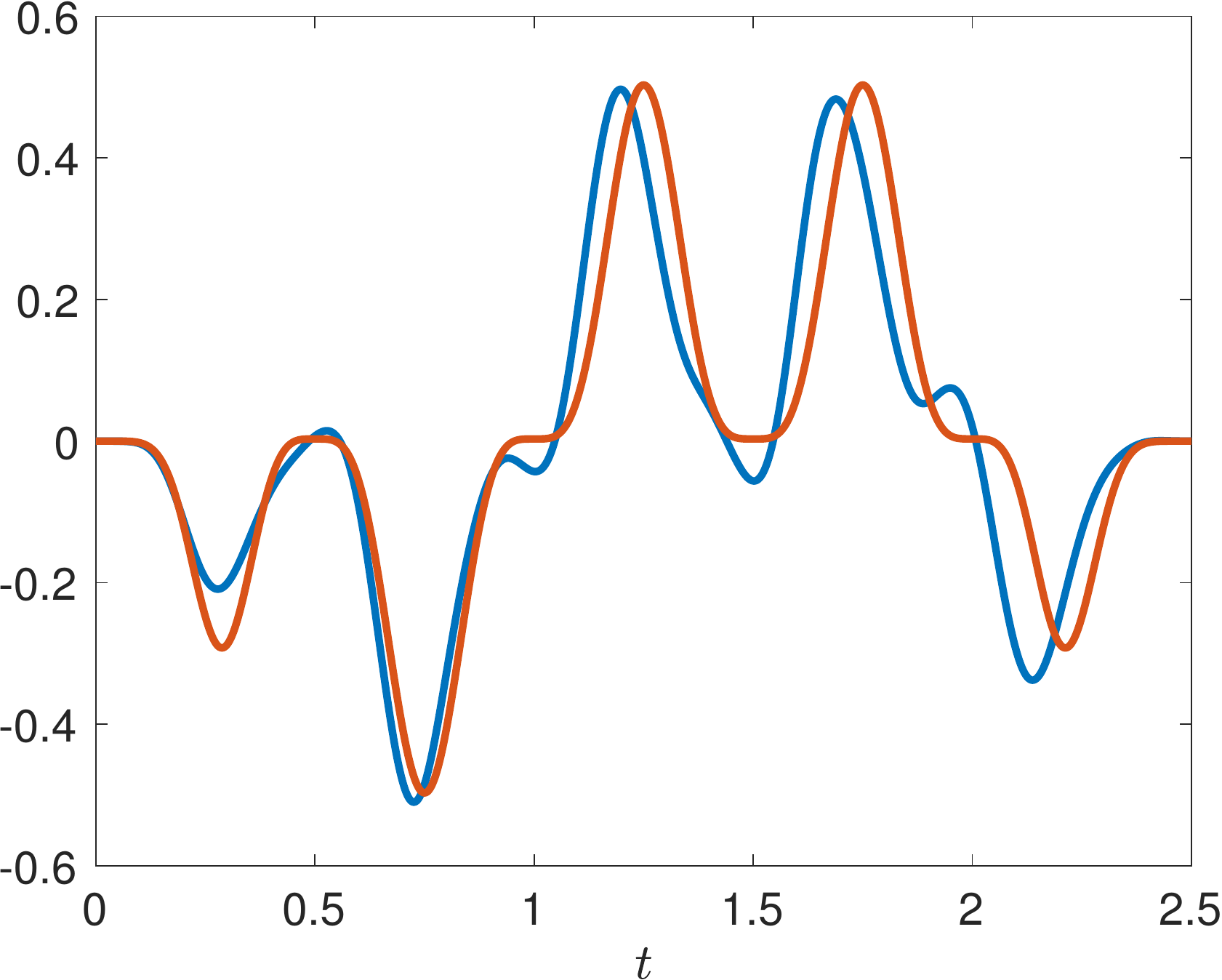}
\caption{Controls $\sqrt{\eps}v^\eps$ (blue) and $v^0$ (red) over $[0,T]$ and $\eps=10^{-1}$ (top-left), $\eps=10^{-2}$ (top-right), $\eps=10^{-3}$ (bottom-left) and $\eps=10^{-4}$ (bottom-right).}\label{fig:veps_v0}
\end{center}
\end{figure}

Figure \ref{fig:veps_v0} depicts the function $\sqrt{\eps}v^{\eps}$ (in blue) on $[0,T]$ for $\eps\in \{10^{-1}10^{-2},10^{-3},10^{-4}\}$	and highlights the punctual convergence of $\sqrt{\eps}v^\eps$ toward $v^0$, the Dirichlet control (in red) for the wave equation and the initial condition $(-y_0,-y_1)$. It is also interesting to note the influence of the amplitude of the parameter $\eps$ on the structure of the control: for $\eps$ large, the control $v^{\eps}$  presents much more oscillations than for $\eps$ small. This property is in agreement with \cite{cindea_beam} where controls for the beam equation are computed.

\begin{figure}[http]
\begin{center}
\includegraphics[scale=0.4]{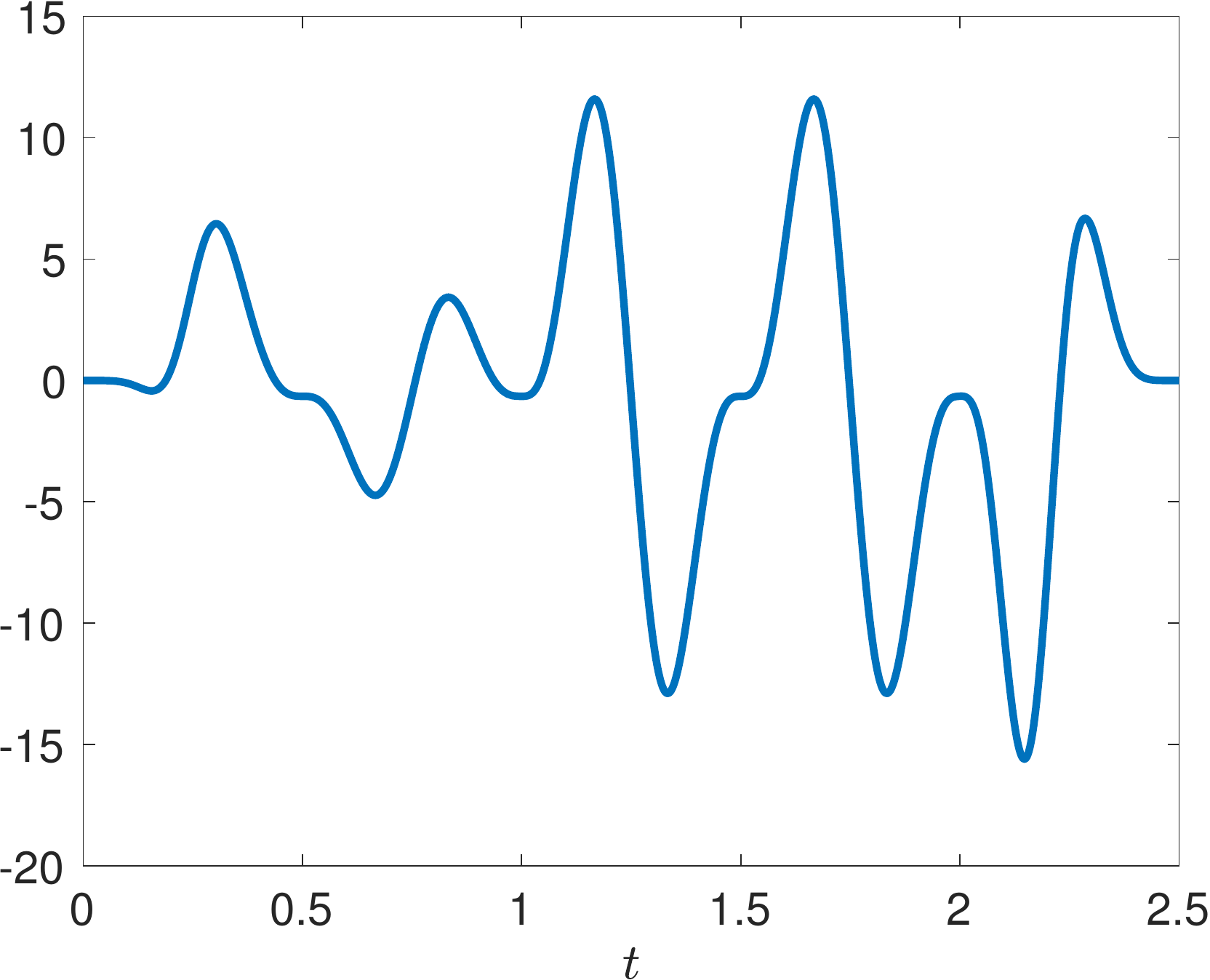}
\hspace*{0.3cm}
\includegraphics[scale=0.4]{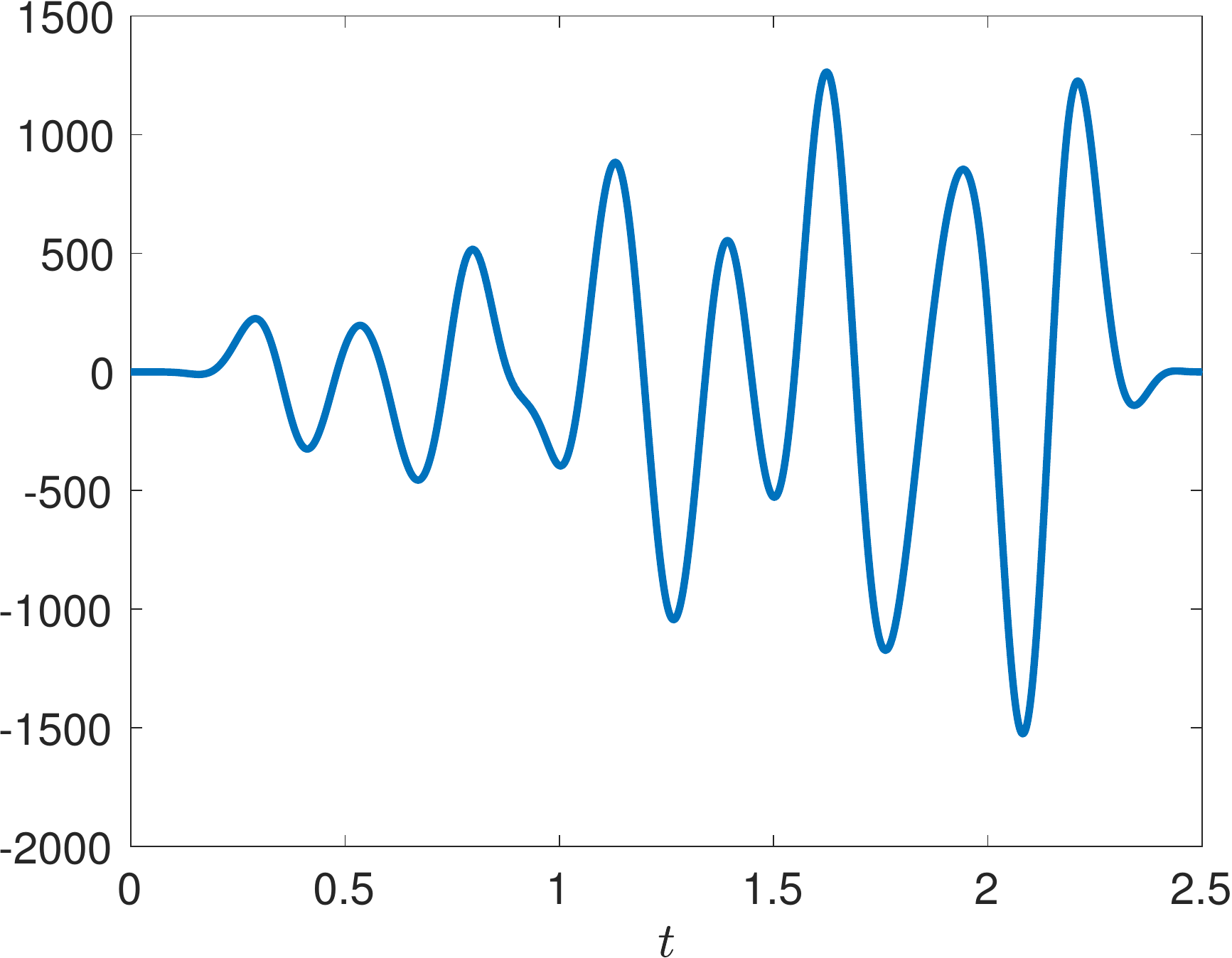}\\
\caption{Controls $v^1$ (left) and $v^2$ (right) over $[0,T]$.}\label{fig:v1_v2}
\end{center}
\end{figure}

\section{Conclusions - Perspectives}

We have rigorously derived an asymptotic expansion of a null control for a singular linear partial differential equation involving a small parameter $\eps>0$. Precisely, we have shown that the control of minimal $L^2$-norm $v^{\eps}$ can be expanded as follows : 
\begin{equation}
v^{\eps}=\frac{1}{\sqrt{\eps}}(v^0+\sqrt{\eps}v^1+\eps v^2) + \mathcal{O}(\eps^{3/4}) \label{devcontrole}
\end{equation}
for the $L^2(0,T)$-norm where the functions $v^k, k=0,1,2$ are related to Dirichlet controls for the wave equation. This strong convergent results requires regularity on the initial data to be controlled, namely $(y_0,y_1)$ in a subset of $H^6(0,1)\times H^4(0,1)$ and refines earlier weak convergence type results given in \cite{LionsPerturbation87}. In particular, we recover that the Neumann control $v^{\eps}$ is singular and that $\sqrt{\eps}v^{\eps}$ converges to a Dirichlet control for the wave equation. It is also important to observe that this singular behavior is not related to the spectral properties of the underlying operator but to the boundary layer occurring on the solution as $\eps$ tends to zero. In particular, a distributed control in the domain does not share a priori such property. To our knowledge, this kind of analysis mixing asymptotic expansion and exact controllability for singular partial differential equation is original. 

From this analysis, a natural question consists to determine the behavior with respect to $\eps$ of the cost of control associated to (\ref{eq:sys_yeps})  and defined as follows 

\begin{equation}  \nonumber 
K(\eps,T):=\sup\limits_{\|y_0,y_1\|_{Z^6}= 1}
\left\{\min_{v\in\mathcal{C}(y_0,y_1,\eps,T)}\| \eta^{-1} v\|_{L^2(0,T)}\right\}
\end{equation}
where $\mathcal{C}(y_0,y_1,\eps,T):=\left\{v\in L^2(0,T); y^\eps=y^\eps(v^\eps)\,\, \textrm{solves}\,\, (\ref{eq:sys_yeps})\,\, \textrm{and satisfies}\,\, (\ref{eq:null}) \right\}$
denotes the non empty set of null controls. In particular, we are looking for the minimal time of controllability defined as $T^\star:=\inf\{T>0; \sup_{\eps>0} K(\eps,T)<\infty\}$
for which the cost is uniformly bounded. The determination of $T^\star$, larger or equal to $2$ in view of Theorem \ref{theoremLions}, is a delicate issue since the initial condition $(y_0,y_1)$ achieving $K(\eps,T)$ may depend on $\eps$. Again, we refer to \cite{CoronGuerrero} which exhibits non intuitive phenomena in the similar context of an advection-diffusion equation. 

The asymptotic expansion of the exact control $v^{\eps}$ is also relevant from an approximation viewpoint, since the expansion \eqref{devcontrole} involves controls for wave equations which are simpler to approximate than $v^\eps$, a fortiori for small values of $\eps$. Moreover, it allows to obtain a convergent approximation of the function $\sqrt{\eps}v^{\eps}$. 
Assume that for $k=0,1,2$, $\{v^k_h\}_{(h>0)}$, is approximation of $v^k$, $h$ being a discretization parameter, satisfying $\Vert v^k-v_h^k\Vert_{L^2(0,T)}=\mathcal{O}(h)$. Such uniform approximation may be achieved using the variational method developed in \cite{cindea_munch_calcolo2015} (see also \cite{cindea_munch_COCV}). Then, in view of (\ref{devcontrole}), the approximation $v^\eps_h:=\eps^{-1/2}(v^0_h+\sqrt{\eps}v_h^1+\eps v_h^2)$ satisfies the estimate 
$$
\Vert \sqrt{\eps}(v^\eps-v^\eps_h)\Vert_{L^2(0,T)}= \mathcal{O}(\eps^{3/4})+\mathcal{O}(h), \quad \forall \eps>0, h>0.
$$

We also mention that the method of matched asymptotic expansions is general and can be used for many other controllability problems involving a small parameter. We mention the model of an elastic cylindric arch (considered in \cite{munch_arch2010}) of length one and constant curvature $c>0$

\begin{equation}
\label{arch_essentialspectrum}
\left\{
\begin{aligned}
& u^{\eps}_{tt}- (u^{\eps}_{x}+c v^{\eps})_x=0, &  \text{ in } Q_T, \\
& v^{\eps}_{tt}+ c (u^{\eps}_{x}+c v^{\eps})+ \eps v^{\eps}_{xxxx}=0, &   \text{ in } Q_T, \\
& u^{\eps}(0,\cdot)=v^{\eps}(0,\cdot)=v^{\eps}_x(0,\cdot)=v^\eps(1,\cdot)=0, &   \text{ in } (0,T), \\
& u^{\eps}(1,\cdot)=f^\eps, v^{\eps}_x(1,\cdot)=g^\eps &   \text{ in } (0,T), \\
& (u^{\eps}(\cdot,0),u^{\eps}_t(\cdot,0))=(u_0,u_1), \, (v^{\eps}(\cdot,0),v^{\eps}_t(\cdot,0))=(v_0,v_1), &  \text{ in } (0,1).
\end{aligned}
\right.
\end{equation}

$u^{\eps}$ and $v^{\eps}$ denote respectively the tangential and normal displacement of the arch. For any $T>0$ large enough, $\eps>0$ and initial conditions $(u_0,u_1)\in H^1_0(0,1)\times L^2(0,1)$, $(v_0,v_1)\in H^2_0(0,1)\times L^2(0,1)$, this system is null controllable through the controls $f^{\eps}$ and $g^{\eps}$. 
The second equation of this system and (\ref{eq:sys_yeps}) share a similar structure: therefore, $v^{\eps}$ exhibits a boundary layer which makes the control $g^{\eps}$ not uniformly bounded with respect to $\eps$. In addition, and contrary to (\ref{eq:sys_yeps}), the underlying operator involves an essential spectrum (as $\eps\to 0$) computed in \cite{geymonat} so that \eqref{arch_essentialspectrum} is not uniformly controllable with respect to the data, as $\eps\to 0$. Nevertheless, we may use the approach developed in this work, and assuming regularity on the data, determine an asymptotic expansion of the two controls $f^{\eps}$ and $g^{\eps}$. 

There are many other partial differential equation involving a small (singular) parameter. We mention the case of the dissipative wave equation  ($\omega$ denotes an open nonempty subset of $(0,1)$)
\begin{equation}\nonumber
\left\{
\begin{aligned}
& \eps y^\eps_{tt}+y^\eps_t-y^\eps_{xx}=v^{\eps}1_{\omega}, & \text{ in } Q_T,\\
& y^\eps(0,t)=y^\eps(1,t)=0, & \text{ in } (0,T), \\
& (y^\eps(\cdot,0),y^\eps_t(\cdot,0))=(y_0,y_1), &\text{ in } (0,1)
\end{aligned}
\right.
\end{equation}
controllable for any $\eps>0$ and for which one can find a sequence of controls $\{v^{\eps}\}_{\eps>0}$ which converges to a null control for the heat equation (we refer to $\cite{lopez_zuazua}$).

\appendix
\section{Appendice}

In this section we prove the convergence results stated in Propositions \ref{pr_1}, \ref{pr_2} and \ref{pr_2}. The proofs require two lemmas that we state first. 

\begin{lemma} \label{le_nhom}
Let $\psi^\varepsilon$ be the solution of the system
\begin{equation} \label{eq_nhom}
\left\{
\begin{array}{ll}
\psi_{tt}^\varepsilon + \varepsilon  \psi^\varepsilon_{xxxx} - \psi^\varepsilon_{xx} =f, & \mbox{in}\,\, Q_T, \\
\psi^\varepsilon (0,\cdot)=g_1, \quad \psi^\varepsilon (1,\cdot)=g_2, & \mbox{in}\quad (0,T),\\
\psi^\varepsilon_{x} (0,\cdot)= h_1, \quad \psi^\varepsilon_{x} (1,\cdot)=h_2, & \mbox{in}\quad  (0,T),\\
\psi^\varepsilon(\cdot,0)=\psi_0, \quad \psi^\varepsilon_t(\cdot,0)=\psi_1, & \mbox{in}\quad \Omega,
\end{array}
\right.
\end{equation}
where $f\in L^1(0,T;L^2)$, $g_1,g_2,h_1,h_2\in H^2(0,T)$ and $(\psi_0,\psi_1)\in H^2\times L^2$ satisfying the compatibility conditions 
\begin{equation} \label{eq_comp1}
\psi^\varepsilon_0(0)=g_1(0), \quad  \psi^\varepsilon_0(1)=g_2(0), \quad  \psi^\varepsilon_{0,x}(0)=h_1(0), \quad  \psi^\varepsilon_{0,x}(1)=h_2(0).
\end{equation}
Then, there exists a constant $C>0$ such that
\begin{equation}
\begin{aligned} 
& \| \psi^\varepsilon \|_{L^\infty(0,T;H^1)} + \| \psi_t^\varepsilon \|_{L^\infty(0,T;L^2)} +\varepsilon^{1/2} \| \psi^\varepsilon_{xx} \|_{L^\infty(0,T;L^2)} \leq C \; F(f,g_1,g_2,h_1,h_2,\psi_0,\psi_1),\\
& \| \varepsilon^{1/2} \psi^\varepsilon_{xx}(0,\cdot) +\psi^\varepsilon_{x} (0,\cdot) \|_{L^2(0,T)} \leq C \; F(f,g_1,g_2,h_1,h_2,\psi_0,\psi_1),\\ \label{eq_est_3}
& \| \varepsilon^{1/2} \psi^\varepsilon_{xx}(1,\cdot) -\psi^\varepsilon_{x} (1,\cdot) \|_{L^2(0,T)} \leq C \; F(f,g_1,g_2,h_1,h_2,\psi_0,\psi_1),
\end{aligned}
\end{equation}
where 
\begin{equation} \nonumber
\begin{aligned}
F(f,g_1,g_2,h_1,h_2,\psi_0,\psi_1)&=\| f\|_{L^1(0,T;L^2)} + \| g_1 \|_{H^2(0,T)}+\| g_2 \|_{H^2(0,T)}\\
&+ \varepsilon^{1/2} (\| h_1 \|_{H^2(0,T)}+\| h_2 \|_{H^2(0,T)})  + \| (\psi_0,\psi_1)\|_{H^1\times L^2}+\varepsilon^{1/2} \| \psi_{0,xx} \|_{L^2}.
\end{aligned}
\end{equation}
\end{lemma}

\textsc{Proof of Lemma} \ref{le_nhom}-
We first homogenize the boundary conditions. Consider the function 
\begin{eqnarray*}
\zeta(x,t)&=&p_1(x)g_1(t) +p_2(x)h_1(t)+ p_3(x)g_2(t) +p_4(x) h_2(t),
\end{eqnarray*}
where $p_i(x)$, $i=1,...,4$ denotes four degree polynomials satisfying
\begin{eqnarray*}
&& p_1(0)=1,\quad p_1'(0)=p_1(1)=p_1'(1)=0,\\
&& p_2'(0)=1,\quad p_2(0)=p_2(1)=p_2'(1)=0,\\
&& p_3(1)=1,\quad p_3(0)=p_3'(0)=p_3'(1)=0,\\
&& p_4'(1)=1,\quad p_4(0)=p_4'(0)=p_4(1)=0.
\end{eqnarray*}

Clearly, $\tilde \psi^\varepsilon(x,t)=\psi^\varepsilon(x,t)-\zeta(x,t)$ satisfies the homogeneous system
\begin{equation} \label{eq_nhom_p}
\left\{
\begin{array}{ll}
\tilde \psi_{tt}^\varepsilon + \varepsilon \tilde \psi^\varepsilon_{xxxx} - \tilde \psi^\varepsilon_{xx} =\tilde f(x,t), & \mbox{in} \quad Q_T, \\
\tilde \psi^\varepsilon (0,t)= \tilde \psi^\varepsilon (1,t)=0, & t\in(0,T)\\
\tilde \psi^\varepsilon_{x} (0,t)=\tilde \psi^\varepsilon_{x} (1,t)=0, & t\in(0,T)\\
\tilde \psi^\varepsilon(x,0)=\tilde \psi_0, \quad \tilde \psi^\varepsilon_t(x,0)=\tilde \psi_1, & x\in\Omega,
\end{array}
\right.
\end{equation}
with 
\begin{eqnarray*}
\tilde f(x,t)&=&f(x,t)-g(x,t)\in L^1(0,T;L^2(0,1)),\\
\tilde \psi_0 &=& \psi_0 -\zeta(x,0)\in H^2_0(0,1), \\
\tilde \psi_1 &=& \psi_1 -\zeta_t(x,0)\in L^2(0,1),
\end{eqnarray*}
where $g = -\zeta_{tt}+\varepsilon\zeta_{xxxx}-\zeta_{xx}$. Note that, 
$$
\|g\|_{L^1(0,T;L^2(0,1))}\leq C(\|g_1^{\prime\prime}\|_{L^1(0,T)}+\|g_2^{\prime\prime}\|_{L^1(0,T)}+\|h_1^{\prime\prime}\|_{L^1(0,T)}+\|h_2^{\prime\prime}\|_{L^1(0,T)}),
$$
for some constant $C$ independent of $\eps$. Then we consider $\tilde \psi^\varepsilon=\tilde \psi^{\varepsilon,1}+\tilde \psi^{\varepsilon,2}$ where $\tilde \psi^{\varepsilon,1}$ is the solution of (\ref{eq_nhom_p}) with $\tilde f(x,t)=0$ and $\tilde \psi^{\varepsilon,2}$ the one associated to $\psi_0=\psi_1=0$. 
Classical energy estimates lead to 
\begin{eqnarray*}
&& \| \tilde \psi^{\varepsilon,1} \|_{L^\infty(0,T;H^1_0)} + \| \tilde \psi_t^{\varepsilon,1} \|_{L^\infty(0,T;L^2)} +\varepsilon^{1/2} \| \tilde \psi^{\varepsilon,1}_{xx} \|_{L^\infty(0,T;L^2)} \\
&& \quad \leq  \| (\psi_0^T,\psi_1^T)\|_{H^1_0\times L^2}+\varepsilon^{1/2} \| \psi^T_{0,xx} \|_{L^2}.
\end{eqnarray*}
On the other hand, by Duhamels formula we obtain 
\begin{eqnarray*}
&& \| \tilde \psi^{\varepsilon,2} \|_{L^\infty(0,T;H^1_0)} + \| \tilde \psi_t^{\varepsilon,2} \|_{L^\infty(0,T;L^2)} +\varepsilon^{1/2} \| \tilde \psi^{\varepsilon,2}_{xx} \|_{L^\infty(0,T;L^2)} \\
&& \quad \leq  C\| \tilde f \|_{L^1(0,T;L^2)}\leq C F(f,g_1,g_2,h_1,h_2,\psi_0,\psi_1).
\end{eqnarray*}
Now, taking into account that 
$\psi^\varepsilon(x,t) = \zeta^\varepsilon(x,t)+\tilde \psi^{\varepsilon,1}+\tilde \psi^{\varepsilon,2},
$ and the above estimates we easily find the first estimate in (\ref{eq_est_3}). The last two estimates in (\ref{eq_est_3}) can be reduced to the corresponding ones for the homogeneous system (\ref{eq_nhom_p}), which is deduced in \cite{LionsPerturbation86}. $\hfill\Box$

\bigskip

\begin{lemma} \label{le_nhom_2}
Let $\psi^\varepsilon$ be the solution of the system
\begin{equation} \label{eq_nhom_2}
\left\{
\begin{array}{ll}
\psi_{tt}^\varepsilon + \varepsilon  \psi^\varepsilon_{xxxx} - \psi^\varepsilon_{xx} =f_1(t) e^{-x/\varepsilon^{1/2}}+f_2(t)e^{-(1-x)/\varepsilon^{1/2}}, & \mbox{in} \,\, Q_T, \\
\psi^\varepsilon (0,\cdot)=\psi^\varepsilon (1,\cdot)=0, & \mbox{in}\quad (0,T),\\
\psi^\varepsilon_{x} (0,\cdot)=\psi^\varepsilon_{x} (1,\cdot)=0, & \mbox{in}\quad  (0,T),\\
\psi^\varepsilon(\cdot,0)=\psi^\varepsilon_t(\cdot,0)=0, & \mbox{in}\quad \Omega,
\end{array}
\right.
\end{equation}
where $f_i\in H^1(0,T)$, $i=1,2$.
Then, there exists a constant $C>0$ such that
\begin{equation} \label{eq_est_1_2}
\begin{aligned}
& \| \psi^\varepsilon \|_{L^\infty(0,T;H^1)} + \| \psi_t^\varepsilon \|_{L^\infty(0,T;L^2)} +\varepsilon^{1/2} \| \psi^\varepsilon_{xx} \|_{L^\infty(0,T;L^2)} \leq C \; \varepsilon^{3/4} F(f_1,f_2),\\ 
& 
\| \varepsilon^{1/2} \psi^\varepsilon_{xx}(0,\cdot) \|_{L^2(0,T)} \leq C \; \varepsilon^{3/4} F(f_1,f_2),\\
&
\| \varepsilon^{1/2} \psi^\varepsilon_{xx}(1,\cdot) \|_{L^2(0,T)} \leq C \; \varepsilon^{3/4} F(f_1,f_2),
\end{aligned}
\end{equation}
where $F(f_1,f_2)=\| f_1\|_{H^1(0,T)} + \| f_2 \|_{H^1(0,T)}$.
\end{lemma}

\textsc{Proof of Lemma} \ref{le_nhom_2}- Without loss of generality we consider the case $f_2=0$ and write $f(t)=f_1(t)$ to simplify. Note that if we apply directly Lemma \ref{le_nhom} to the solution of (\ref{eq_nhom_2}) we simply obtain estimates for the left hand side in (\ref{eq_est_1_2}) that depend on the $L^1(0,T;L^2(0,1))$-norm of $f(t)e^{-x/\varepsilon^{1/2}}$, i.e.
$$
\| f(t)e^{-x/\varepsilon^{1/2}} \|_{L^1(0,T;L^2(0,1))} =\|f\|_{L^1(0,T)} \| e^{-x/\varepsilon^{1/2}}\|_{L^2(0,1)} \leq C \|f\|_{L^1(0,T)} \varepsilon^{1/4}.
$$
Here we use an energy argument to improve these estimates up to the power $\varepsilon^{3/4}$.

Multiplying the first equation in system (\ref{eq_nhom_2}) by $\psi^\varepsilon_t$ and integrating in space, we obtain 
\begin{eqnarray*}
\frac{dE^\varepsilon (t)}{dt} &=& f(t)\int_0^1  \psi^\varepsilon_t e^{-x/\varepsilon^{1/2}} dx=\frac{d}{dt} \left[ f(t)\int_0^1  \psi^\varepsilon e^{-x/\varepsilon^{1/2}} dx \right]-f'(t)\int_0^1  \psi^\varepsilon e^{-x/\varepsilon^{1/2}} dx\\
&=&\frac{d}{dt} \left[ f(t)\varepsilon^{1/2}\int_0^1  \psi^\varepsilon_x e^{-x/\varepsilon^{1/2}} dx \right]-f'(t)\varepsilon^{1/2}\int_0^1  \psi^\varepsilon_x e^{-x/\varepsilon^{1/2}} dx, 
\end{eqnarray*}
where 
$$
E^\varepsilon(t)= \int_0^1 \left[|\psi^\varepsilon_t|^2+\varepsilon |\psi^\varepsilon_{xx}|^2+|\psi^\varepsilon_{x}|^2 \right]dx.
$$
Integrating now in time we obtain, for all $t\in (0,T)$, 
\begin{eqnarray} \nonumber
E^\varepsilon(t)&=&\varepsilon^{1/2}f(t)\int_0^1  \psi^\varepsilon_x e^{-x/\varepsilon^{1/2}} dx-\varepsilon^{1/2}\int_0^t f'(s)\int_0^1  \psi^\varepsilon_x(x,s) e^{-x/\varepsilon^{1/2}} dx ds \\ \label{eq_le_me}
&\leq& \varepsilon^{1/2} |f(t)| (E^{\varepsilon}(t))^{1/2} \| e^{-x/\varepsilon^{1/2}}\|_{L^2} + \varepsilon^{1/2} \| e^{-x/\varepsilon^{1/2}}\|_{L^2} \| f' \|_{L^2(0,t)} \left(\int_0^t E(s)ds \right)^{1/2}.
\end{eqnarray}
and therefore, 
$$
\int_0^T E^\varepsilon(t) dt \leq \varepsilon^{1/2} \| e^{-x/\varepsilon^{1/2}}\|_{L^2(0,1)} \big(\|f\|_{L^2(0,T)}+T\| f' \|_{L^2(0,T)}\big)\left(\int_0^T E(t)dt \right)^{1/2}.
$$
Therefore,
$$
\int_0^T E^\varepsilon(t) dt \leq \varepsilon \| e^{-x/\varepsilon^{1/2}}\|_{L^2(0,1)}^2 (\|f\|_{L^2(0,T)}+T\| f' \|_{L^2(0,T)})^2.
$$
Substituting in (\ref{eq_le_me}) and taking into account that $\| e^{-x/\varepsilon^{1/2}}\|_{L^2(0,1)}\leq C \varepsilon^{1/4}$ we obtain
$$
\sup_{t\in[0,T]} E^{\varepsilon}(t)\leq C \varepsilon^{3/4}\| f'\|_{L^2} \sup_{t\in[0,T]} (E^{\varepsilon}(t))^{1/2} + \varepsilon^{3/2} \| f^{\prime}\|_{L^2(0,T)}(\|f\|_{L^2(0,T)}+T\| f' \|_{L^2(0,T)}),
$$
and we easily deduce (\ref{eq_est_1_2}). 

We now prove the last two estimates in (\ref{eq_est_1_2}). We only consider the first one since the other is analogous. Multiplying the first equation in system (\ref{eq_nhom_2}) by $x\psi^\varepsilon_x$ and integrating we easily obtain 
\begin{eqnarray*}
&&  \int_0^T\!\!\!\int_0^1 f(t) e^{-x/\varepsilon^{1/2}} x\psi^\varepsilon_x \; dxdt = \left. \int_0^1 \psi^\varepsilon_t x \psi^\varepsilon_x dx \right]_0^T \\
&& \quad +\int_0^T\!\!\!\int_0^1  \left( -\psi^\varepsilon_t x \psi^\varepsilon_{xt}-\varepsilon \psi^\varepsilon_{xxx} (\psi^\varepsilon_{x}+x \psi^\varepsilon_{xx}) + \frac{x}{2}  (|\psi^\varepsilon_{x}|^2)_x \right) dxdt\\
&& = \int_0^1 \psi^\varepsilon_t(x,T) x \psi^\varepsilon_x(x,T) dx + \int_0^T\!\!\!\int_0^1 \left( |\psi^\varepsilon_t|^2+ \frac32 \varepsilon |\psi^\varepsilon_{xx}|^2 - \frac12
|\psi^\varepsilon_{x}|^2 \right) dxdt \\
&& \quad - \left. \int_0^T \frac12 x |\psi^\varepsilon_{xx}|^2 dt \right]_{x=0}^{x=1} .
\end{eqnarray*}
Therefore, we have
\begin{eqnarray*}
&& \frac12 \int_0^T |\psi^\varepsilon_{xx}(1,t)|^2 dt =\int_0^1 \psi^\varepsilon_t(x,T) x \psi^\varepsilon_x(x,T) dx \\
&& \quad+ \int_0^T\!\!\!\int_0^1 \left( |\psi^\varepsilon_t|^2+ \frac32 \varepsilon |\psi^\varepsilon_{xx}|^2 - \frac12
|\psi^\varepsilon_{x}|^2 -f(t) e^{-x/\varepsilon^{1/2}} x\psi^\varepsilon_x\right) dxdt \\
&& \leq \| \psi^\varepsilon_t(\cdot,T) \|_{L^2(0,1)} \| \psi^\varepsilon_x(\cdot,T) \|_{L^2(0,1)} + \frac32 \int_0^T E^\varepsilon(t) dt + \|f\|_{L^2(0,T)} \| e^{-x/\varepsilon^{1/2}}\|_{L^2(0,1)} \| \psi^\varepsilon_x \|_{L^\infty([0,T];L^2(0,1))}\\
&& \leq C_1 \| E^\varepsilon \|_{L^\infty([0,T])} + C_2 \|f\|_{L^2(0,T)} \varepsilon^{1/4} \| (E^\varepsilon)^{1/2} \|_{L^\infty([0,T])}. 
\end{eqnarray*}

\begin{equation}
\begin{aligned}
 \frac12 \int_0^T |\psi^\varepsilon_{xx}(1,t)|^2 dt &=\int_0^1 \psi^\varepsilon_t(x,T) x \psi^\varepsilon_x(x,T) dx \\
& \qquad+ \int_0^T\!\!\!\int_0^1 \left( |\psi^\varepsilon_t|^2+ \frac32 \varepsilon |\psi^\varepsilon_{xx}|^2 - \frac12
|\psi^\varepsilon_{x}|^2 -f(t) e^{-x/\varepsilon^{1/2}} x\psi^\varepsilon_x\right) dxdt \\
& \leq \| \psi^\varepsilon_t(\cdot,T) \|_{L^2(0,1)} \| \psi^\varepsilon_x(\cdot,T) \|_{L^2(0,1)} + \frac32 \int_0^T E^\varepsilon(t) dt \\
& \qquad + \|f\|_{L^2(0,T)} \| e^{-x/\varepsilon^{1/2}}\|_{L^2(0,1)} \| \psi^\varepsilon_x \|_{L^\infty([0,T];L^2(0,1))}\\
& \leq C_1 \| E^\varepsilon \|_{L^\infty([0,T])} + C_2 \|f\|_{L^2(0,T)} \varepsilon^{1/4} \| (E^\varepsilon)^{1/2} \|_{L^\infty([0,T])}. 
\end{aligned}
\end{equation}

The result then follows from this estimate and (\ref{eq_est_1_2}). $\hfill\Box$
 
\bigskip

\par\noindent
{\sc Proof of Proposition \ref{pr_1}.} We introduce $\psi^\varepsilon= y^\varepsilon-y^{\varepsilon,0}$ solution of (\ref{eq_nhom}) 
with 
\begin{equation} \label{eq_le31_7}
\begin{aligned}
& f=-v^0_{tt}(t)e^{-(1-x)/\varepsilon^{1/2}}-\varepsilon y^0_{xxxx}, \quad g_1= -v^0 e^{-1/\varepsilon^{1/2}}, \, g_2=0,\\ 
& h_1= -y^0_x(0,\cdot)-\eps^{-1/2}e^{-1/\varepsilon^{1/2}} v^0 , \quad h_2= v^\varepsilon-\varepsilon^{-1/2}v^0-y^0_x(1,\cdot), \\
& (\psi_0^\varepsilon,\psi_1^\varepsilon)= (y_0^{\varepsilon}, y_1^{\varepsilon})-(y_0^{0}, y_1^{0}).
\end{aligned}
\end{equation}
By the hypothesis on the regularity of both the initial and boundary data, the solution $y^0$ of system (\ref{eq_sis0}) satisfies 
$$
y^0\in C([0,T]; H^4(0,1))\cap  C^1([0,T]; H^3(0,1)).
$$
In particular, 
$$
y^0_{xxxx}\in L^1(0,T;L^2(0,1)), \quad y^0_x(0,\cdot), \; y^0_x(1,\cdot) \in H^2(0,T).
$$
The compatibility conditions (\ref{eq_comp1}) are also easily verified. For example, 
\begin{eqnarray*}
&& \psi^\varepsilon_0(0)=y_0^\varepsilon(0)-y_0(0)=0=-v^0(0)e^{-1/\varepsilon^{1/2}}=g_1(0), \\
&& \psi^\varepsilon_{0,x}(0)=y_{0,x}^\varepsilon(0)-y_{0,x}(0)=-y_{0,x}(0) =-y_{0,x}(0)-v^0(0)\varepsilon^{-1/2} e^{-1/\varepsilon^{1/2}}=h_1(0),
\end{eqnarray*}
and similarly for those at $x=1$. 
Therefore the result is direct consequence of Lemmas \ref{le_nhom} and \ref{le_nhom_2}. $\hfill\Box$

\bigskip

We now prove Proposition \ref{pr_2}. The proof of Propsition \ref{pr_3} is analogous and we omit it. 

\bigskip

\par\noindent
{\sc Proof of Proposition \ref{pr_2}-} Following the idea in the proof of Proposition \ref{pr_2} we try to apply Lemma \ref{le_nhom} to $\psi^\varepsilon=y^\varepsilon - y^{\varepsilon,1}$. We obtain that $\psi^\varepsilon$ is solution of  (\ref{eq_nhom}) with
 \begin{equation}
 \begin{aligned}
&f=-\varepsilon y^0_{xxxx}(x,t)+ y^0_{tt}(1,t) e^{-(1-x)/\sqrt{\varepsilon}} -\varepsilon^{3/2} y^1_{xxxx}(x,t)  \\
& \hspace*{2cm} + \varepsilon^{1/2} y^1_{tt}(0,t) e^{-x/\varepsilon^{1/2}} +\varepsilon^{1/2} y^1_{tt}(1,t) e^{-(1-x)/\sqrt{\varepsilon}}, \\ 
&  g_1=(y^0(1,\cdot)+\varepsilon^{1/2}y^1(1,\cdot))e^{-1/\varepsilon^{1/2}}, \quad g_2=(y^0(0,\cdot)+\varepsilon^{1/2}y^1(0,\cdot))e^{-1/\varepsilon^{1/2}},\\ 
& h_1=-\varepsilon^{1/2}y^1_x(0,\cdot)+(\varepsilon^{-1/2}y^0(1,\cdot)+y^1(1,\cdot))e^{-1/\varepsilon^{1/2}},\\
& h_2= v^{\eps}-\eps^{-1/2}(v^0+\sqrt{\eps}v^1)-\varepsilon^{1/2}y^1_x(1,\cdot)-y^1(0,\cdot)e^{-1/\eps^{1/2}}, \\
& (\psi_0,\psi_1)= {\cal O}(\varepsilon).
\end{aligned}
\end{equation} 
By the regularity hypothesis on the initial and boundary data, and the compatibility conditions in (\ref{eq_compo1}), the solutions $y^0$ and $y^1$ of systems (\ref{eq_sis0}) and (\ref{eq_sis11}) satisfy 
$$
y^j\in C([0,T]; H^{5-j}(0,1))\cap  C^1([0,T]; H^{4-j}(0,1)),\quad j=0,1.
$$
We now write $\psi^\varepsilon=\psi^\varepsilon_1+\psi^\varepsilon_2$, where $\psi^\varepsilon_1$ satisfies (\ref{eq_nhom}) with $f=-\varepsilon y^0_{xxxx}-\varepsilon^{3/2} y^1_{xxxx}$ and the same boundary and initial conditions as $\psi^\varepsilon$. Then $\psi^\varepsilon_2$ will be solution of (\ref{eq_nhom_2}) with second hand term given by 
$$
(y^0_{tt}(1,t)+\varepsilon^{1/2} y^1_{tt}(1,t)) e^{-(1-x)/\sqrt{\varepsilon}} + \varepsilon^{1/2} y^1_{tt}(0,t) e^{-x/\varepsilon^{1/2}}. 
$$ 

Thus, we can apply Lemmas \ref{le_nhom} and \ref{le_nhom_2} to $\psi^\varepsilon_1$ and $\psi^\varepsilon_2$ respectively. Combining both estimates allow to conclude the proof.  
$\hfill\Box$

\bibliographystyle{siam}


\end{document}